\documentclass[envcountsect]{svjour3} 

\smartqed
\usepackage{graphicx}
\usepackage{amsmath}
\usepackage{amsfonts}
\usepackage{amssymb}
\usepackage{bm}
\journalname{JOTA}

\begin{document}

\title{The Game of Two Identical Cars: an Analytical Description of the Barrier}

\author{Maksim Buzikov \and Andrey Galyaev }

\institute{Maksim Buzikov, Corresponding author \at
             V. A. Trapeznikov Institute of Control Sciences of Russian Academy of Sciences \\
              Moscow, Russia\\
              me.buzikov@physics.msu.ru
           \and
           Andrey Galyaev \at
              V. A. Trapeznikov Institute of Control Sciences of Russian Academy of Sciences \\
              Moscow, Russia\\
              galaev@ipu.ru
}

\date{Received: date / Accepted: date}

\maketitle

\begin{abstract}
    In this paper, a pursuit-evasion game of two players known as the game of two identical cars is examined. It is assumed that the game proceeds on a two-dimensional plane. Both players have constant speeds and limited turn radii. The goal of the first player (pursuer) is to ensure that the second player (evader) is guaranteed to be within a capture circle as quickly as possible. The goal of the evader is to avoid the capture or delay it as long as possible. The kinematics of both players are described by the same equations. Thus, the game has only one free parameter, the capture radius. This work aims at a fully analytical description of the barrier surface for all values of the capture radius. Previously, A.W.~Merz analytically investigated the barrier of the game of two identical cars. In this work, it is found that there is a certain critical value of the capture radius, above which the barrier occurs qualitatively differently from Merz's example. Also, we obtained an analytical description of the optimal feedback controls on the barrier.
\end{abstract}

\section{Introduction}

The game of two cars (GTC) was originally described and investigated by R.~Isaacs in the 1950s. In his seminal work \cite[pp.~237--244]{isaacs1965differential}, Isaacs obtained a partial description of the barrier for this problem and a description of the singular lines on the barrier. Isaacs supposed that the independent parameters of the problem (capture radius, speed ratio, and minimum-turn radius ratio) are such that the capture is possible from an arbitrary initial state. Thus, the barrier explored by Isaacs is not closed. In \cite{cockayne1967plane}, E.~Cockayne revealed necessary and sufficient conditions on the parameters of the problem of interception from an arbitrary initial state. These conditions were refined by G.T.~Rublein in \cite{rublein1972pursuit}. Necessary and sufficient conditions for the existence of an evasion strategy are presented in \cite{borowko1984game}. The publications \cite{cockayne1967plane,rublein1972pursuit,borowko1984game} suppose that interception is a coincide of the pursuer and evader planar coordinates. It corresponds to a zero capture radius. In \cite{breakwell1977minimum}, J.V.~Breakwell and A.W.~Merz calculated minimal capture radius, which always permits the interception from an arbitrary initial state for the given speed ratio and minimum-turn radius ratio. A new criterion of full state-space capturability based on the geometry of the barrier was proposed in \cite{pachter1980geometry}. It states that if the origin of the coordinate system is bonded with the pursuer (the y-axis is aligned with the velocity of the pursuer) and the barrier surface doesn't cross the positive y-axis then the pursuer can intercept the evader from any initial state. Later, M.~Pachter and T.~Miloh revealed new regions of the parameters space where the open barrier surface qualitatively differs from previously observed kinds \cite{pachter1987geometric}. An optimal feedback control synthesis for the GTC was introduced in \cite{simakova1967differential,simakova1983optimal} by E.N.~Simakova. Also, Simakova investigated the value of the game in \cite{simakova1986bellman} and the game of kind in \cite{simakova1970problem}. Her results are applicable only in a special case when the pursuer surpasses all the parameters of the evader. Also, she considers only a state-space region where the value of the game is a smooth function. In \cite{kostsov1988differential}, Simakova examines obtained results in comparison with other methods such as proportional navigation. In \cite{farber1980approximate}, N.~Farber and J.~Shinar described a methodology to obtain an approximate solution and compared it with Simakova's solution. A computational approach based on discretization of the GTC is presented in \cite[pp.~48--60]{baron1970new}. In \cite{bera2017comprehensive}, R.~Bera et al. investigated and visualized singular surfaces of the game in a wide range of parameters. A comprehensive survey of the GTC in the context of pursuit-evasion differential games is presented in \cite{weintraub2020introduction}.

The solution of the GTC can be utilized for a collision avoidance problem where it is supposed that there is no cooperation between the players and the evader has to imply a strategy such that its use will lead to escape for any strategy of the pursuer. Usually, the parameters of the collision avoidance problem differ from a pursuit-evasion problem in that the pursuer doesn't have the advantage of the speed and minimum-turn radius. In \cite{meier1969new}, L.~Meier introduced a method that helps to obtain a capture region when it is bounded. His method is based on an analysis of the reachable sets. Initially, the pursuer's and evader's reachable sets in a given time are determined. Further, by combining these loci for various fixed times and a fixed initial state and finding envelopes and intersections, the capture region is determined. In publications \cite{merz1973optimal,ciletti1976collision,merz1976collision}, A.~Merz et al. examined a cooperative collision avoidance problem with equations of motion of the GTC. Cooperative collision avoidance implies maximization of distance to the closest point of approach by both players instead maximization of time to collision by the evader in the pursuit-evasion game. Further, T.~Miloh and S.D.~Sharma supplemented Merz's results and investigated the barrier for a collision avoidance problem for the cooperative maneuvers of players \cite{miloh1975determination,sharma1977ship}. In \cite{vincent1972problem,vincent1974some}, T.L.~Vincent et al. analyzed the problem of the collision avoidance of two vehicles when there is no cooperation of players. The main idea of the analysis is to divide the state space into some zones of vulnerability. The "red" zone corresponds to a win region of the pursuer in the GTC. So, a boundary of the red zone is the barrier of the GTC. In \cite{vincent1972problem}, the analytical description of the barrier in case of evader's superiority in parameters is provided. In \cite{vincent1976aircraft}, the same authors supplemented their analysis by considering a proportional navigation strategy of the pursuer. Papers \cite{olsder1978differential,miloh1983game} examine a collision-avoidance problem for line segment and elliptical shape players. In \cite{miloh1976maritime,miloh1989ship}, the extended kinematic model taking into account variable speeds of players is considered. In \cite{kwik1989calculation}, K.H.~Kwik calculated an escape maneuver for the evader if the pursuer maneuver is known.

The game of two identical cars (GTIC) was originally considered by A.~Merz in \cite{merz1972game}. In the GTIC, it is assumed that the speeds and minimum-turn radii of the players are the same. After choosing an appropriate description of the GTIC statement there is only one parameter in the problem, the capture radius. In \cite{merz1972game}, A.~Merz described the barrier of the game, the universal and dispersal surfaces. A game of surveillance evasion of two identical cars is considered by I.~Greenfeld in \cite{greenfeld1987differential}. It differs from the GTIC in the game-set, it is inside of the capture circle, and also it differs in that the pursuer maximizes while the evader minimizes the time to escape. In \cite{mitchell2001games}, I.~Mitchell uses the description of the barrier obtained by Merz to validate a solver of the Hamilton-Jacobi equation. A cooperative collision avoidance problem of two identical cars is completely analyzed by T.~Tarnopolskaya and N.~Fulton in \cite{tarnopolskaya2009optimal}. Further, these authors examined more general cases of unequal parameters of players in \cite{tarnopolskaya2010synthesis,tarnopolskaya2010dispersal,tarnopolskaya2010non,tarnopolskaya2011synthesis,maurer2012optimal,maurer2012singular,maurer2014computation}.

Although the GTIC has been studied intensively, there are some unexplored questions regarding this problem at the moment. The analytical description of the barrier presented in \cite{merz1972game,mitchell2001games} is constrained only by the parametric equations of surfaces that constitute the barrier. These publications do not provide analytical conditions on the limitations of the parameter domains in the description of surfaces, i.e., there are no conditions for cutting off excess parts of these surfaces. Similar applies to the description of the universal and dispersal surfaces from \cite{merz1972game}. In addition, Merz noted an unexpected result that occurs in the numerical calculation of the angular slice of the dispersal surface at a small enough angle between the velocities of the players. Analytical analysis of this phenomenon is also not covered in the available sources.

This paper is devoted to a complete analytical description of the barrier. The analysis shows that for different values of the capture radius, there is a different analytical description of the conditions for cutting off the excess parts of the surfaces that form the barrier. Instead of a trial and error approach from \cite{mitchell2001games}, we present a systematic way of constructing the barrier. This paper also proposes efficient methods for calculating optimal feedback controls for players and provides an analytical expression for the barrier parametrization in terms of the state vector.

The rest of the paper is organized as follows. Section~\ref{sec:problem_formulation} provides the mathematical statement of the GTIC. Section~\ref{sec:adjoint_strategies} describes the strategies of players in terms of the state vector and the normal vector. Next, Section~\ref{sec:barrier} is devoted to an analytical description of the barrier and its form depending on the capture radius. Then Section~\ref{sec:optimal_feedback_controls} provides a synthesis of optimal feedback controls for players and proposes an analytical description of the barrier in terms of the state vector. Finally, Section~\ref{sec:conclusions} gives conclusions and future work.

\section{Problem formulation}\label{sec:problem_formulation}

In this section, we formulate the mathematical statement of the GTIC. We use the notation $\boldsymbol{z} = \begin{bmatrix} x & y & \theta \end{bmatrix}^\mathrm{T} \in \mathbb{R}^2 \times \mathbb{S}$ for a configuration, which is a planar position and orientation triplet. The space $\mathbb{S}$ consists of the real numbers, but the equality in $\mathbb{S}$ differs from the equality in $\mathbb{R}$. For every $\theta_1, \theta_2 \in \mathbb{S}$, equality $\theta_1 = \theta_2$ holds in $\mathbb{S}$ when there exists a number $k \in \mathbb{Z}$ such that equality $\theta_1 = \theta_2 + 2\pi k$ holds in $\mathbb{R}$. Throughout this paper, the subscripts $P$ and $E$ refer to the pursuer and evader configuration correspondingly.

\begin{figure}
    \begin{center}
        \includegraphics[width=0.4\linewidth]{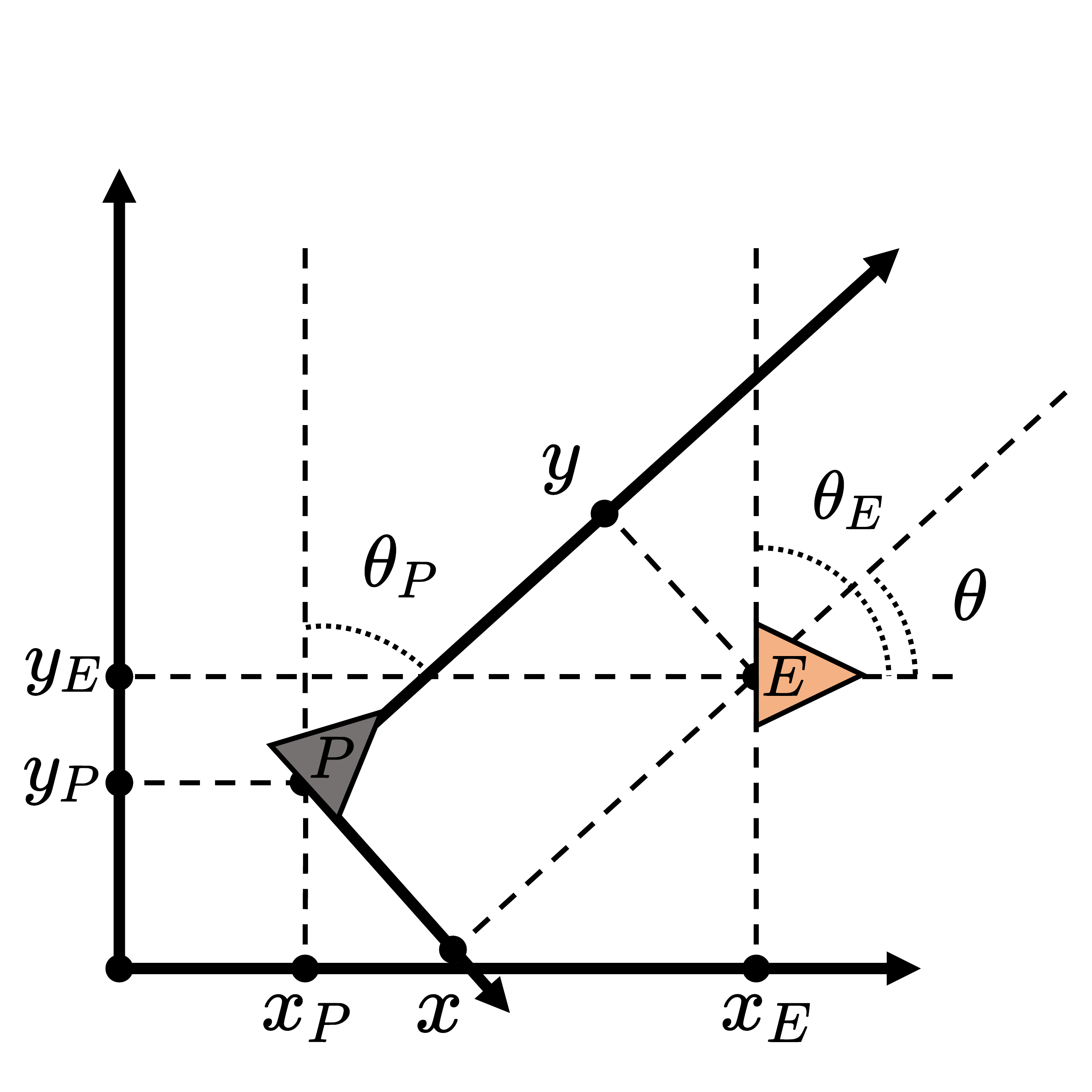}
        \caption{Players in the realistic space $\begin{bmatrix} x_P & y_P & \theta_P & x_E & y_E & \theta_E\end{bmatrix}^\mathrm{T}$ and reduced space $\begin{bmatrix} x & y & \theta \end{bmatrix}^\mathrm{T}$}
        \label{fig:coordinate_system}
    \end{center}
\end{figure}

For the GTIC, the kinematics of each player is described by the Dubins model. According to this model, pursuer's and evader's configurations in realistic space could be written as $\begin{bmatrix} x_P & y_P & \theta_P \end{bmatrix}^\mathrm{T}$ and $\begin{bmatrix} x_E & y_E & \theta_E \end{bmatrix}^\mathrm{T}$, respectively (see Fig.~\ref{fig:coordinate_system}).
All angles count in a clockwise direction of the y-axis. In reduced space, we associate the coordinate system with the pursuer's position and align the y-axis along with the velocity of the pursuer. The equations of motion of the GTIC~\cite[Eq. 2]{merz1972game} are presented as
\begin{equation}\label{eq:dynamics}
    \dot{\boldsymbol{z}} 
       = 
        \begin{bmatrix}
               \dot{x} \\
               \dot{y} \\
               \dot{\theta}
        \end{bmatrix}
        =
        \begin{bmatrix}
               -uy + \sin\theta\\
               -1 + ux + \cos\theta\\
               v - u
        \end{bmatrix}
    \overset{\mathrm{def}}{=} \boldsymbol{f}(\boldsymbol{z}, \boldsymbol{u}, \boldsymbol{v}).
\end{equation}
Here, the pursuer's and evader's controls correspond to functions
\begin{equation*}
    u(\boldsymbol{z}) \overset{\mathrm{def}}{=} \boldsymbol{u}(\boldsymbol{z}) \in \mathcal{U} \overset{\mathrm{def}}{=} [-1, +1], \quad  v(\boldsymbol{z}) \overset{\mathrm{def}}{=} \boldsymbol{v}(\boldsymbol{z})\in \mathcal{V} \overset{\mathrm{def}}{=} [-1, +1],
\end{equation*}
where $\mathcal{U}$ and $\mathcal{V}$ are the constrained ranges of players' controls.

The set of all possible states defined as the game set is given by
\begin{equation*}
    \mathcal{S} \overset{\mathrm{def}}{=} \{\boldsymbol{z} \in \mathbb{R}^2 \times \mathbb{S}:\: x^2 + y^2 \geq \ell^2\},
\end{equation*}
where $\ell \in \mathbb{R}^+$ is a capture radius of the target set
\begin{equation*}\label{terminal_set}
    \mathcal{C} \overset{\mathrm{def}}{=} \partial\mathcal{S} = \{\boldsymbol{z} \in \mathbb{R}^2 \times \mathbb{S}:\: x^2 + y^2 = \ell^2\}.
\end{equation*}
We suppose that the game starts from the initial state $\boldsymbol{z}_0 \in \mathcal{S}$ at the time moment $t = 0$, then the current state of the game $\boldsymbol{z}(t) = \begin{bmatrix} x(t) & y(t) & \theta(t)\end{bmatrix}^\mathrm{T}$ depends on the players' controls and it can be defined as a solution of the equations of motion \eqref{eq:dynamics}. If the state reaches the target set $\mathcal{C}$ for the first time at some moment $t_f \in \mathbb{R}^+_0$, then the game terminates. Further, we associate the subscript $f$ with some final states throughout this paper. If the target set cannot be reached for the given controls, then we suppose that $t_f = +\infty$.

The terminal state $\boldsymbol{z}(t_f)$ on the target set can be represented in a more convenient form. Let $\boldsymbol{s} = \begin{bmatrix} \varphi_f & \theta_f \end{bmatrix}^\mathrm{T} \in \mathbb{S}^2$. The target set can be presented as
\begin{equation*}
    \mathcal{C} = \{\boldsymbol{z}_\mathcal{C}(\boldsymbol{s}):\: \boldsymbol{s}\in \mathbb{S}^2\}, \text{ where } \boldsymbol{z}_\mathcal{C}(\boldsymbol{s}) \overset{\mathrm{def}}{=} \begin{bmatrix} \ell\sin\varphi_f & \ell\cos\varphi_f & \theta_f \end{bmatrix}^\mathrm{T}_.
\end{equation*}
The polar angle $\varphi_f$ fixes the planar position on the capture circle.

The outcome functional of the game is
\begin{equation*}
    J[\boldsymbol{u}(\cdot), \boldsymbol{v}(\cdot); \boldsymbol{z}_0] \overset{\mathrm{def}}{=} \int\limits_0^{t_f}L(\boldsymbol{z}(t), \boldsymbol{u}(\boldsymbol{z}(t)), \boldsymbol{v}(\boldsymbol{z}(t)))\mathrm{d}t + G(\boldsymbol{z}(t_f)),
\end{equation*}
where the running cost function $L(\boldsymbol{z}, \boldsymbol{u}, \boldsymbol{v}) = 1$, and the terminal cost function $G(\boldsymbol{z}) = 0$. The pursuer strives to minimize the outcome functional when the evader maximizes it. For the optimal feedback controls $\boldsymbol{u}^*(\cdot)$, $\boldsymbol{v}^*(\cdot)$ and for any other admissible controls $\boldsymbol{u}(\cdot)$, $\boldsymbol{v}(\cdot)$ the following inequalities hold
\begin{equation*}
    J[\boldsymbol{u}^*(\cdot), \boldsymbol{v}(\cdot); \boldsymbol{z}_0] \leq J[\boldsymbol{u}^*(\cdot), \boldsymbol{v}^*(\cdot); \boldsymbol{z}_0] \leq J[\boldsymbol{u}(\cdot), \boldsymbol{v}^*(\cdot); \boldsymbol{z}_0].
\end{equation*}
The value function of the game is defined by
\begin{equation*}
    V(\boldsymbol{z}_0) \overset{\mathrm{def}}{=} J[\boldsymbol{u}^*(\cdot), \boldsymbol{v}^*(\cdot); \boldsymbol{z}_0].
\end{equation*}
The value $V(\boldsymbol{z}_0)$ is a time to capture if the pursuer and the evader play optimally.

It is obvious that there is a set of states where the evader can escape with any strategy of the pursuer. For example, if the velocities of players are collinear ($\theta = 0_\mathbb{S}$)\footnote{We will use the symbol $0_\mathbb{S}$ to denote a zero element of space $\mathbb{S}$.} and the evader is in front of the pursuer, then the control $v(\boldsymbol{z}) = 0$ guarantees evasion. It follows, that the escape set $\mathcal{F} \subset \mathcal{S}$ is not empty. We denote by $\mathcal{W} \subset \mathcal{S}$ a capture set. The boundary that separates the capture set from the escape set is a barrier $\mathcal{B} \subset \mathcal{S}$. In the rest of the paper, we describe the barrier $\mathcal{B}$ analytically. Also, we provide the controls $\boldsymbol{u}^*(\boldsymbol{z})$, $\boldsymbol{v}^*(\boldsymbol{z})$ for $\boldsymbol{z} \in \mathcal{B}$.

\section{Main equation}\label{sec:adjoint_strategies}

If an initial state of the game lies on barrier $\mathcal{B}$ and both players use the optimal strategies, then the state vector doesn't leave the barrier. It means that the evader must be able to set a value of its control in such a way that it prevents the state vector from penetrating the capture set $\mathcal{W}$ while the pursuer keeps the state out of the escape set $\mathcal{F}$ in the same way. The described property to prevent penetration constitutes the notion of semi-permeability.

Let $\nu_{\boldsymbol{z}} = \begin{bmatrix} \nu_x & \nu_y & \nu_\theta \end{bmatrix}^\mathrm{T} \in \mathbb{R}^3$ denote a normal vector of barrier $\mathcal{B}$ at some point. The length of $\nu_{\boldsymbol{z}}$ is not necessarily a unit. According to Isaacs \cite[pp.~205--210]{isaacs1965differential}, the condition that $\mathcal{B}$ be a semipermeable surface
\begin{equation}\label{eq:main_equation}
    \min_{\boldsymbol{u} \in \mathcal{U}} \max_{\boldsymbol{v} \in \mathcal{V}}\:(\nu_{\boldsymbol{z}}, \boldsymbol{f}(\boldsymbol{z}, \boldsymbol{u}, \boldsymbol{v})) = 0
\end{equation}
for all $\boldsymbol{z} \in \mathcal{B}$. Using the equations of motion~\eqref{eq:dynamics} gives 
\begin{equation*}
    (\nu_{\boldsymbol{z}}, \boldsymbol{f}(\boldsymbol{z}, \boldsymbol{u}, \boldsymbol{v})) = -u(y\nu_x - x\nu_y + \nu_\theta) + v\nu_\theta + \nu_x\sin\theta + \nu_y\cos\theta - \nu_y.
\end{equation*}
In the particular case of the GTIC, the right part of equations of motion is separable and it can be presented in a special form
\begin{equation*}
    \boldsymbol{f}(\boldsymbol{z}, \boldsymbol{u}, \boldsymbol{v}) = \boldsymbol{f}_P(\boldsymbol{z}, \boldsymbol{u}) + \boldsymbol{f}_E(\boldsymbol{z}, \boldsymbol{v}) + \boldsymbol{g}(\boldsymbol{z}),
\end{equation*}  
where we have set
\begin{equation*}
     \boldsymbol{f}_P(\boldsymbol{z}, \boldsymbol{u}) \overset{\mathrm{def}}{=} \begin{bmatrix} -uy \\ ux \\ -u \end{bmatrix}, \quad \boldsymbol{f}_E(\boldsymbol{z}, \boldsymbol{v}) \overset{\mathrm{def}}{=} \begin{bmatrix} 0 \\ 0 \\ v \end{bmatrix}, \quad \boldsymbol{g}(\boldsymbol{z}) \overset{\mathrm{def}}{=} \begin{bmatrix} \sin\theta \\ -1 + \cos\theta \\ 0 \end{bmatrix}.
\end{equation*}
The functions $\boldsymbol{f}_P$, $\boldsymbol{f}_E$ are linear in components of the controls for both players. It follows that the main equation~\eqref{eq:main_equation} can be analyzed through the investigation of switch functions \cite[pp.~147--148]{lewin2012differential}. We will denote by
\begin{equation*}
    s_P(\boldsymbol{z}, \nu_{\boldsymbol{z}}) \overset{\mathrm{def}}{=} y\nu_x - x\nu_y + \nu_\theta, \quad s_E(\boldsymbol{z}, \nu_{\boldsymbol{z}}) \overset{\mathrm{def}}{=} \nu_\theta
\end{equation*}
the pursuer's and evader's switch functions. The main equation~\eqref{eq:main_equation} helps to determine players' optimal controls depending on the state $\boldsymbol{z}$ and the normal vector $\nu_{\boldsymbol{z}}$
\begin{equation}\label{eq:implicit_controls}
    \begin{split}
        &\tilde{u}(\boldsymbol{z},\nu_{\boldsymbol{z}}) \overset{\mathrm{def}}{=} \mathrm{arg}\min_{\boldsymbol{u} \in \mathcal{U}}\:(\nu_{\boldsymbol{z}}, \boldsymbol{f}_P(\boldsymbol{z}, \boldsymbol{u})),\\
        &\tilde{v}(\boldsymbol{z},\nu_{\boldsymbol{z}}) \overset{\mathrm{def}}{=} \mathrm{arg}\max_{\boldsymbol{v} \in \mathcal{V}}\:(\nu_{\boldsymbol{z}}, \boldsymbol{f}_E(\boldsymbol{z}, \boldsymbol{v})).
    \end{split}
\end{equation}
Applying the notions of switch functions we obtain 
\begin{equation}\label{eq:switch_functions_suppl}
    (\boldsymbol{\nu}_{\boldsymbol{z}}, \boldsymbol{f}_P(\boldsymbol{z}, \boldsymbol{u})) = -us_P(\boldsymbol{z},\nu_{\boldsymbol{z}}), \quad (\boldsymbol{\nu}_{\boldsymbol{z}}, \boldsymbol{f}_E(\boldsymbol{z}, \boldsymbol{v})) = vs_E(\boldsymbol{z},\nu_{\boldsymbol{z}}).
\end{equation}

The definitions~\eqref{eq:implicit_controls} can be ambiguous if the corresponding switch function vanishes. In this case, we use the property of admissible controls. It supposes that any admissible control is a piecewise continuous function of the time. Since every piecewise continuous function has a left-continuous equivalent with respect to Lebesgue's measure, we can consider only left-continuous functions. Thus, we demand that $\tilde{u}(\boldsymbol{z}(\cdot), \nu_{\boldsymbol{z}}(\cdot))$ and $\tilde{v}(\boldsymbol{z}(\cdot), \nu_{\boldsymbol{z}}(\cdot))$ are the left-continuous functions of time:
\begin{equation}\label{eq:controls_left_continuity}
    \begin{split}
        &\tilde{u}(\boldsymbol{z}(t), \nu_{\boldsymbol{z}}(t)) = \lim_{\tau \to +0} \tilde{u}(\boldsymbol{z}(t-\tau), \nu_{\boldsymbol{z}}(t-\tau)),\\
        &\tilde{v}(\boldsymbol{z}(t), \nu_{\boldsymbol{z}}(t)) = \lim_{\tau \to +0} \tilde{v}(\boldsymbol{z}(t-\tau), \nu_{\boldsymbol{z}}(t-\tau)).
    \end{split}
\end{equation}
The construction of semipermeable surface $\mathcal{B}$ uses the evolution of normal vector $\nu_{\boldsymbol{z}}$ \cite[pp.~113--114]{lewin2012differential} described by the adjoint equations
\begin{equation}\label{eq:adjoint_equations}
    \dot{\nu}_{\boldsymbol{z}}
        \overset{\mathrm{def}}{=}
        \begin{bmatrix}
               \dot{\nu}_x \\
               \dot{\nu}_y \\
               \dot{\nu}_\theta
        \end{bmatrix}
        = -\frac{\partial\boldsymbol{f}(\boldsymbol{z}, \boldsymbol{u}, \boldsymbol{v})}{\partial \boldsymbol{z}}\nu_{\boldsymbol{z}}
        =
        \begin{bmatrix}
            -u\nu_y\\
            u\nu_x\\
            -\nu_x\cos\theta + \nu_y\sin\theta
        \end{bmatrix}.
\end{equation}
These equations hold along barrier paths. Let $\nu_{\boldsymbol{z}}(t) = \begin{bmatrix} \nu_x(t) & \nu_y(t) & \nu_\theta(t) \end{bmatrix}^\mathrm{T}$ denote a solution of this system at the time moment $t$.

The definitions~\eqref{eq:implicit_controls} and conditions~\eqref{eq:controls_left_continuity}--\eqref{eq:adjoint_equations} help to obtain explicit forms of $\tilde{u}(\cdot)$, $\tilde{v}(\cdot)$.
\begin{lemma}\label{lem:candidate_optimal_controls}
    In the GTIC, the controls $\tilde{u}(\cdot)$, $\tilde{v}(\cdot)$ given by
    \begin{equation*}
        \tilde{u}(\boldsymbol{z}, \nu_{\boldsymbol{z}}) = 
        \begin{cases}
        	&\mathrm{sgn}\: s_P(\boldsymbol{z}, \nu_{\boldsymbol{z}}), \quad s_P(\boldsymbol{z}, \nu_{\boldsymbol{z}}) \neq 0;\\
        	&\mathrm{sgn}\:\nu_x, \quad s_P(\boldsymbol{z}, \nu_{\boldsymbol{z}}) = 0,\: \nu_x \neq 0;\\
        	&0, \quad s_P(\boldsymbol{z}, \nu_{\boldsymbol{z}}) = 0, \: \nu_x = 0, \: \nu_y < 0;\\
        	&u \in \{-1, 0, +1\}, \quad \text{otherwise},\\
        \end{cases}
    \end{equation*}
    \begin{equation*}
        \tilde{v}(\boldsymbol{z}, \nu_{\boldsymbol{z}}) = 
        \begin{cases}
        	&\mathrm{sgn}\:s_E(\boldsymbol{z}, \nu_{\boldsymbol{z}}), \quad s_E(\boldsymbol{z}, \nu_{\boldsymbol{z}}) \neq 0;\\
        	&\mathrm{sgn}(\nu_x\cos\theta - \nu_y\sin\theta), \quad s_E(\boldsymbol{z}, \nu_{\boldsymbol{z}}) = 0, \: \nu_y\sin\theta \neq \nu_x\cos\theta;\\
        	&0, \quad s_E(\boldsymbol{z}, \nu_{\boldsymbol{z}}) = 0, \: \nu_y\sin\theta = \nu_x\cos\theta, \: \nu_x\sin\theta < -\nu_y\cos\theta;\\
        	&v \in \{-1, 0, +1\}, \quad \text{otherwise}, 
        \end{cases}
    \end{equation*}
    satisfy the necessary conditions of optimality.
\end{lemma}
\begin{proof}
    Since $\boldsymbol{z}(\cdot)$, $\nu_{\boldsymbol{z}}(\cdot)$ are differentiable functions of time, the switch function $s_P(\boldsymbol{z}(\cdot), \nu_{\boldsymbol{z}}(\cdot))$ is differentiable too. Using~\eqref{eq:dynamics} and~\eqref{eq:adjoint_equations}, we may compute derivatives of the switch function
    \begin{equation*}
        \dot{s}_P(\boldsymbol{z}(t), \nu_{\boldsymbol{z}}(t)) = -\nu_x(t), \quad \ddot{s}_P(\boldsymbol{z}(t), \nu_{\boldsymbol{z}}(t)) = \tilde{u}(\boldsymbol{z}(t), \nu_{\boldsymbol{z}}(t))\nu_y(t).
    \end{equation*}
    Consider $s_P(\boldsymbol{z}(t), \nu_{\boldsymbol{z}}(t)) \neq 0$. Substituting~\eqref{eq:switch_functions_suppl} into~\eqref{eq:implicit_controls} we rewrite~\eqref{eq:implicit_controls} as
    \begin{equation*}
        \tilde{u}(\boldsymbol{z}(t), \nu_{\boldsymbol{z}}(t)) = \mathrm{sgn}\:s_P(\boldsymbol{z}(t), \nu_{\boldsymbol{z}}(t)) = \mathrm{sgn}(y(t)\nu_x(t) - x(t)\nu_y(t) + \nu_\theta(t)).
    \end{equation*}
    We now turn to the case $s_P(\boldsymbol{z}(t), \nu_{\boldsymbol{z}}(t)) = 0$. According to the mean value theorem, for all $\tau \geq 0$, there exists $\alpha \in (0, 1)$ such that
    \begin{equation*}
        s_P(\boldsymbol{z}(t-\tau), \nu_{\boldsymbol{z}}(t-\tau)) = -\tau\dot{s}_P(\boldsymbol{z}(t - \alpha\tau), \nu_{\boldsymbol{z}}(t - \alpha\tau)) = \tau \nu_x(t - \alpha\tau).
    \end{equation*}
    Next, if $\nu_{x}(t) \neq 0$, then there exists a neighborhood where $\nu_x(\cdot)$ remains its sign. Applying~\eqref{eq:controls_left_continuity} we deduce that
    \begin{equation*}
        \tilde{u}(\boldsymbol{z}(t), \nu_{\boldsymbol{z}}(t)) = \lim_{\tau \to +0} \mathrm{sgn}( \tau\nu_x(t - \alpha\tau)) = \mathrm{sgn}(\nu_x(t)).
    \end{equation*}
    Now, set $\nu_x(t) = 0$. According to Lagrange's version of Taylor's theorem, for all $\tau \geq 0$, there exists $\alpha \in (0, 1)$ such that
    \begin{multline*}
        s_P(\boldsymbol{z}(t-\tau), \nu_{\boldsymbol{z}}(t-\tau)) = \frac{\tau^2}{2}\ddot{s}_P(\boldsymbol{z}(t - \alpha\tau), \nu_{\boldsymbol{z}}(t - \alpha\tau))\\
        = \frac{\tau^2}{2}\tilde{u}(\boldsymbol{z}(t - \alpha\tau), \nu_{\boldsymbol{z}}(t - \alpha\tau))\nu_y(t - \alpha\tau).
    \end{multline*}
    Applying~\eqref{eq:controls_left_continuity} yields
    \begin{equation*}
        \tilde{u}(\boldsymbol{z}(t), \nu_{\boldsymbol{z}}(t)) = \lim_{\tau \to +0} \mathrm{sgn}\: \tilde{u}(\boldsymbol{z}(t-\alpha\tau), \nu_{\boldsymbol{z}}(t-\alpha\tau))\nu_y(t - \alpha\tau).
    \end{equation*}
    If $\nu_y(t) < 0$, then there exists a neighborhood where $\nu_y(\cdot)$ remains negative. In this case we see that 
    \begin{equation*}
        \tilde{u}(\boldsymbol{z}(t), \nu_{\boldsymbol{z}}(t)) = -\lim_{\tau \to +0} \mathrm{sgn}\: \tilde{u}(\boldsymbol{z}(t-\alpha\tau), \nu_{\boldsymbol{z}}(t-\alpha\tau)) = 0
    \end{equation*}
    The same conclusion can be drawn for $\nu_y(t) > 0$. There exists a neighborhood where $\nu_y(\cdot)$ remains positive. It implies
    \begin{equation*}
        \tilde{u}(\boldsymbol{z}(t), \nu_{\boldsymbol{z}}(t)) = \lim_{\tau \to +0} \mathrm{sgn}\: \tilde{u}(\boldsymbol{z}(t-\alpha\tau), \nu_{\boldsymbol{z}}(t-\alpha\tau)) \in \{-1, 0, +1\},
    \end{equation*}
    since the range of $\mathrm{sgn}(\cdot)$ is $\{-1, 0, +1\}$. 
    
    The same reasoning applies to the explicit form of $\tilde{v}(\cdot)$. Derivatives of the switch function are equal to 
    \begin{equation*}
        \begin{split}
            &\dot{s}_E(\boldsymbol{z}(t), \nu_{\boldsymbol{z}}(t)) = -\nu_x(t)\cos\theta(t) + \nu_y(t)\sin\theta(t),\\
            &\ddot{s}_E(\boldsymbol{z}(t), \nu_{\boldsymbol{z}}(t)) = \tilde{v}(\boldsymbol{z}(t), \nu_{\boldsymbol{z}}(t))\left(\nu_x(t)\sin\theta(t) + \nu_y(t)\cos\theta(t)\right).
        \end{split}
    \end{equation*}
    We leave it to the reader to obtain the provided explicit form of $\tilde{v}(\cdot)$. \qed
\end{proof}

Obtained explicit forms of controls satisfy the necessary conditions of optimality. Hence, $\tilde{u}(\cdot)$, $\tilde{v}(\cdot)$ set only candidate optimal control laws.

In the remainder of this section, we give solutions of~\eqref{eq:dynamics} and~\eqref{eq:adjoint_equations} for fixed values $u$, $v$ of controls \cite{merz1972game,mitchell2001games}. These solutions are required for the investigation of the barrier's emanation in the next section. It is more convenient to explore the emanation using the retrograde time variable $\tau = t_f - t$. Let $\tilde{\boldsymbol{z}} = \begin{bmatrix} \tilde{x} & \tilde{y} & \tilde{\theta} \end{bmatrix}^\mathrm{T} \in \mathbb{R}^2 \times \mathbb{S}$. Integration of~\eqref{eq:dynamics} using the fixed values of controls $u \in \mathcal{U}$, $v \in \mathcal{V}$ and condition $\boldsymbol{z}(t_f) = \tilde{\boldsymbol{z}}$ leads to the following solution:
\begin{multline}\label{eq:trajectory_for_fixed_controls}
    \boldsymbol{z}_{u, v}(\tau; \tilde{\boldsymbol{z}}) \overset{\mathrm{def}}{=} \boldsymbol{z}(t_f - \tau)\\
    =
    \begin{bmatrix}
        \tilde{x}\cos u\tau + \tilde{y}\sin u\tau + \frac{u\tau^2}2\mathrm{sinc}^2\frac{u\tau}2 - \tau\mathrm{sinc}\frac{v\tau}2\sin\left(\tilde{\theta} + \left(u - \frac{v}2\right)\tau\right)\\
        \tilde{y}\cos u\tau - \tilde{x}\sin u\tau + \tau\mathrm{sinc}\frac{u\tau}2 - \tau\mathrm{sinc}\frac{v\tau}2\cos\left(\tilde{\theta} + \left(u - \frac{v}2\right)\tau\right)\\
        \tilde{\theta} + (u - v)\tau
    \end{bmatrix}.
\end{multline}
The distance between players on the plane is equal to
\begin{equation}\label{eq:radial_distance}
    r_{u, v}(\tau; \tilde{\boldsymbol{z}}) \overset{\mathrm{def}}{=} \sqrt{(\boldsymbol{z}_{u, v}(\tau; \tilde{\boldsymbol{z}}), \boldsymbol{I}_{XY}\boldsymbol{z}_{u, v}(\tau; \tilde{\boldsymbol{z}}))}, \text{ where }
    \boldsymbol{I}_{XY} \overset{\mathrm{def}}{=}
    \begin{bmatrix}
        1 & 0 & 0\\
        0 & 1 & 0\\
        0 & 0 & 0
    \end{bmatrix}.
\end{equation}
Let $\tilde{\nu}_{\boldsymbol{z}} = \begin{bmatrix} \tilde{\nu}_{x} & \tilde{\nu}_{y} & \tilde{\nu}_{\theta} \end{bmatrix}^\mathrm{T} \in \mathbb{R}^3$ and $\nu_{\boldsymbol{z}}(t_f) = \tilde{\nu}_{\boldsymbol{z}}$. By the same way, integrating~\eqref{eq:adjoint_equations} leads to
\begin{multline*}
    \nu_{\boldsymbol{z}}^{u, v}(\tau; \boldsymbol{z}, \nu_{\boldsymbol{z}}) \overset{\mathrm{def}}{=} \nu_{\boldsymbol{z}}(t_f - \tau)\\
    =
    \begin{bmatrix}
        \nu_x\cos u\tau + \nu_y\sin u\tau\\
        \nu_y\cos u\tau - \nu_x\sin u\tau\\
        \nu_{\theta} + \tau\mathrm{sinc}\frac{v\tau}2\left(\nu_x\cos\left(\theta-\frac{v\tau}2\right) - \nu_y\sin\left(\theta - \frac{v\tau}2\right)\right)
    \end{bmatrix}.
\end{multline*}
This solution also requires that the controls are fixed values.

\section{Barrier}\label{sec:barrier}

The first step of the barrier construction is to find the usable part and its boundary. By definition, the usable part is the part of the target set where the pursuer can guarantee termination regardless of the choice of the control $\boldsymbol{v}$ by the evader \cite[pp.~39--40]{lewin2012differential}. The unit normal of the target set $\mathcal{C}$ directed into the game set $\mathcal{S}$ at some point $\boldsymbol{z}_\mathcal{C}(\boldsymbol{s})$ for $\boldsymbol{s} \in \mathbb{S}^2$ will be denoted by $\boldsymbol{n}_\mathcal{C}(\boldsymbol{s})$. According to the definition,
\begin{equation*}
    \left[\frac{\partial\boldsymbol{z}_{\mathcal{C}}(\boldsymbol{s})}{\partial\boldsymbol{s}}\right]^\mathrm{T} \boldsymbol{n}_{\mathcal{C}}(\boldsymbol{s}) =
    \begin{bmatrix}
        \ell\cos\varphi_f & -\ell\sin\varphi_f & 0\\
        0 & 0 & 1
    \end{bmatrix}
    \boldsymbol{n}_{\mathcal{C}}(\boldsymbol{s})
    = \boldsymbol{0}.
\end{equation*}
Therefore, $\boldsymbol{n}_{\mathcal{C}}(\boldsymbol{s}) = \begin{bmatrix}\sin\varphi_f & \cos\varphi_f & 0\end{bmatrix}^\mathrm{T} \in \mathbb{R}^3$. The pursuer can guarantee termination when the velocity vector penetrates to the target set. Hence, 
\begin{multline*}
        \mathcal{UP} \overset{\mathrm{def}}{=} \left\{\boldsymbol{z}_{\mathcal{C}}(\boldsymbol{s}):\: \boldsymbol{s} \in \mathbb{S}^2,\: \min_{\boldsymbol{u} \in \mathcal{U}}\max_{\boldsymbol{v} \in \mathcal{V}} (\boldsymbol{n}_{\mathcal{C}}(\boldsymbol{s}), \boldsymbol{f}(\boldsymbol{z}_{\mathcal{C}}(\boldsymbol{s}), \boldsymbol{u}, \boldsymbol{v})) < 0\right\}\\
        = \left\{\boldsymbol{z}_\mathcal{C}(\boldsymbol{s}):\: \boldsymbol{s} \in \mathbb{S}^2,\: \cos(\varphi_f-\theta_f) - \cos\varphi_f < 0\right\}\\
        =\left\{\boldsymbol{z}_\mathcal{C}(\boldsymbol{s}):\: \theta_f \in (0, 2\pi),\: \varphi_f \in \left(\frac{\theta_f}2 - \pi, \frac{\theta_f}2\right)\right\}.
\end{multline*}

The boundary of usable part ($\mathcal{BUP}$) is a locus of points where the barrier starts its emanation in backward time $\tau = t_f - t$. By definition, $\mathcal{BUP} = \partial\mathcal{UP}$. In the GTIC, the boundary of usable part consist of three lines: $\mathcal{BUP} = \mathcal{BUP}_0 \cup \mathcal{BUP}_{-1} \cup \mathcal{BUP}_{+1}$. The line $\mathcal{BUP}_0$ corresponds to $\theta_f = 0_\mathbb{S}$:
\begin{equation*}
    \mathcal{BUP}_0 \overset{\mathrm{def}}{=} \left\{\boldsymbol{z}_{\mathcal{BUP}_0}(\varphi_f): \varphi_f \in \mathbb{S}\right\} = \{\boldsymbol{z} \in \mathbb{R}^2 \times \mathbb{S}:\: \theta = 0_\mathbb{S}, \: \ell = \sqrt{x^2 + y^2}\},
\end{equation*}
where
\begin{equation*}
    \boldsymbol{z}_{\mathcal{BUP}_0}(\varphi_f) \overset{\mathrm{def}}{=} \boldsymbol{z}_\mathcal{C}(\varphi_f, 0_\mathbb{S}) = \begin{bmatrix} \ell\sin\varphi_f & \ell\cos\varphi_f & 0_{\mathbb{S}} \end{bmatrix}^\mathrm{T}.
\end{equation*}
For future purposes, we specify the unit normal of the target set on $\mathcal{BUP}_0$:
\begin{equation*}
    \boldsymbol{n}_{\mathcal{BUP}_0}(\varphi_f) \overset{\mathrm{def}}{=} \boldsymbol{n}_\mathcal{C}(\varphi_f, 0_{\mathbb{S}}).
\end{equation*}
The last two parts of $\mathcal{BUP}$ correspond to cases $\varphi_f = \theta_f/2 - \pi$ and $\varphi_f = \theta_f/2$. Let us denote\footnote{Throughout the paper, $\mathbb{B}$ stands for the binary set $\{-1, +1\}$. Also, we will use the letter $\upsilon$ (upsilon) to denote a binary variable. This letter looks like evader's control $v$ and it will play this role in fact.} by $\mathcal{BUP}_\upsilon$ these two lines $\upsilon \in \mathbb{B}$. The parametric description of these lines is given by
\begin{equation*}
    \mathcal{BUP}_\upsilon \overset{\mathrm{def}}{=} \left\{\boldsymbol{z}_{\mathcal{BUP}_\upsilon}(\theta_f):\: \theta_f \in (0, 2\pi)\right\},
\end{equation*}
where
\begin{equation*}
    \boldsymbol{z}_{\mathcal{BUP}_\upsilon}(\theta_f) \overset{\mathrm{def}}{=} \boldsymbol{z}_\mathcal{C}\left(\frac{(1 + \upsilon)\pi + \theta_f}2, \theta_f\right) = \begin{bmatrix} -\upsilon\ell\sin\frac{\theta_f}2 & -\upsilon\ell\cos\frac{\theta_f}2 & \theta_f \end{bmatrix}^\mathrm{T}.
\end{equation*}
The unit normal of the target set on $\mathcal{BUP}_\upsilon$ is defined by
\begin{equation*}
    \boldsymbol{n}_{\mathcal{BUP}_\upsilon}(\theta_f) \overset{\mathrm{def}}{=} \boldsymbol{n}_\mathcal{C}\left(\frac{(1 + \upsilon)\pi + \theta_f}2, \theta_f\right).
\end{equation*}

The emanation of the barrier surface is determined by the candidate optimal control laws $\tilde{u}(\cdot)$, $\tilde{v}(\cdot)$ from Lemma~\ref{lem:candidate_optimal_controls}. Since the values of the candidate optimal control laws are ambiguous, all necessary scenarios should be examined. The initial conditions of the retrogressive emanation are given by $\boldsymbol{z}(t_f) = \boldsymbol{z}_\mathcal{C}(\boldsymbol{s}) \in \mathcal{BUP}$, $\nu_{\boldsymbol{z}}(t_f) = \boldsymbol{n}_\mathcal{C}(\boldsymbol{s})$. We consider the emanation from $\mathcal{BUP}$ separately for $\mathcal{BUP}_0$ and $\mathcal{BUP}_\upsilon$.

\subsection{Primary barrier emanation}

For $\mathcal{BUP}_\upsilon$, the candidate optimal control laws $\tilde{u}(\cdot)$ and $\tilde{v}(\cdot)$ are determined certainly and according to Lemma~\ref{lem:candidate_optimal_controls} equal to
\begin{equation}\label{eq:BUPv_controls}
    \begin{aligned}
        &\tilde{u}(\boldsymbol{z}_{\mathcal{BUP}_\upsilon}(\theta_f), \boldsymbol{n}_{\mathcal{BUP}_\upsilon}(\theta_f)) = -\upsilon,\\
        &\tilde{v}(\boldsymbol{z}_{\mathcal{BUP}_\upsilon}(\theta_f), \boldsymbol{n}_{\mathcal{BUP}_\upsilon}(\theta_f)) = \upsilon,
    \end{aligned}
\end{equation}
where $\theta_f \in (0, 2\pi)$. We associate the notion $\mathcal{B}_\mathcal{P}^\upsilon$ with the part of barrier $\mathcal{B}$ which emanates from $\mathcal{BUP}_\upsilon$ (remember that $\upsilon \in \mathbb{B}$). Since not all emanating parts are valid parts of the barrier, we would mark emanating parts with a tilde. Therefore, $\tilde{\mathcal{B}}_\mathcal{P}^\upsilon$ is an emanating from the $\mathcal{BUP}_\upsilon$ surface, and $\mathcal{B}_\mathcal{P}^\upsilon = \tilde{\mathcal{B}}_\mathcal{P}^\upsilon \cap \mathcal{B}$ is a valid part of the barrier. The redundant part $\tilde{\mathcal{B}}_\mathcal{P}^\upsilon \setminus \mathcal{B}_\mathcal{P}^\upsilon$ will be determined later. Combining~\eqref{eq:BUPv_controls} with the solution of equations of motion for fixed values of controls \eqref{eq:trajectory_for_fixed_controls} we obtain
\begin{multline*}
    \boldsymbol{z}_{\mathcal{B}_\mathcal{P}^\upsilon}(\tau, \theta_f) \overset{\mathrm{def}}{=} \boldsymbol{z}_{-\upsilon, \upsilon}(\tau; \boldsymbol{z}_{\mathcal{BUP}_\upsilon}(\theta_f))\\
    =
    \begin{bmatrix}
        -\upsilon\ell\sin\left(\frac{\theta_f}2 - \upsilon\tau\right) - \upsilon\left(1 - \cos\tau + \cos(\theta_f - 2\upsilon\tau) - \cos(\theta_f - \upsilon\tau)\right)\\
        -\upsilon\ell\cos\left(\frac{\theta_f}2 - \upsilon\tau\right) + \sin{\tau} + \upsilon\left(\sin(\theta_f - 2\upsilon\tau) - \sin(\theta_f - \upsilon\tau)\right)\\
        \theta_f - 2\upsilon\tau
    \end{bmatrix}.
\end{multline*}
The normal vector of $\tilde{\mathcal{B}}_\mathcal{P}^\upsilon$ at $\boldsymbol{z}_{\mathcal{B}_\mathcal{P}^\upsilon}(\tau, \theta_f)$ is given by
\begin{equation*}
    \nu_{\boldsymbol{z}, \mathcal{B}_\mathcal{P}^\upsilon}(\tau, \theta_f) \overset{\mathrm{def}}{=} \nu^{-\upsilon, \upsilon}_{\boldsymbol{z}}(\tau; \boldsymbol{z}_{\mathcal{BUP}_\upsilon}(\theta_f), \boldsymbol{n}_{\mathcal{BUP}_\upsilon}(\theta_f)).
\end{equation*}
The range of change for the retrograde time $\tau$ on $\tilde{\mathcal{B}}_\mathcal{P}^\upsilon$ starts from $\tau = 0$ and ends at an arbitrary moment. We can define this time moment in such a way that it is guaranteed not less than either the minimal switch time (when one of the switch functions vanishes) or the time when $\tilde{\mathcal{B}}_\mathcal{P}^\upsilon$ intersects another part of $\mathcal{B}$. The second option is not available now since the other parts of $\mathcal{B}$ are not described. To compute the minimal switch time amounts to finding a minimal positive root of
\begin{align*}
    s_P(\boldsymbol{z}_{\mathcal{B}_\mathcal{P}^\upsilon}(\tau, \theta_f), \nu_{\boldsymbol{z}, \mathcal{B}_\mathcal{P}^\upsilon}(\tau, \theta_f)) = \cos\frac{\theta_f}2 - \cos\left(\frac{\theta_f}2 - \upsilon\tau\right) = 0,\\
    s_E(\boldsymbol{z}_{\mathcal{B}_\mathcal{P}^\upsilon}(\tau, \theta_f), \nu_{\boldsymbol{z}, \mathcal{B}_\mathcal{P}^\upsilon}(\tau, \theta_f)) = \cos\left(\frac{\theta_f}2 - \upsilon\tau\right) - \cos\frac{\theta_f}2 = 0.
\end{align*}
A trivial verification shows that the root is $\tau = (1 - \upsilon)\pi + \upsilon\theta_f$ for $\theta_f \in (0, 2\pi)$. Note that the third component of $\boldsymbol{z}_{\mathcal{B}_\mathcal{P}^\upsilon}(\tau, \theta_f)$ becomes zero when
\begin{equation*}
    \tau = \frac{(1 - \upsilon)\pi + \upsilon\theta_f}2, \quad \theta_f \in (0, 2\pi),
\end{equation*}
which is earlier than the switch functions vanish. Since the evader can escape in the case $\theta = 0_\mathbb{S}$ duplicating a strategy of the pursuer, the corresponding state doesn't belong to $\mathcal{B}$ and lies on the redundant part of $\tilde{\mathcal{B}}_\mathcal{P}^\upsilon$. Summarizing, we set
\begin{equation*}
    \tilde{\mathcal{B}}^\upsilon_\mathcal{P} \overset{\mathrm{def}}{=} \left\{\boldsymbol{z}_{\mathcal{B}_\mathcal{P}^\upsilon}(\tau, \theta_f):\: \theta_f \in (0, 2\pi),\: \tau \in \left(0, \frac{(1 - \upsilon)\pi + \upsilon\theta_f}2\right)\right\}.
\end{equation*}
For future purposes, we use another parametrization of $\tilde{\mathcal{B}}_\mathcal{P}^\upsilon$ given by a new variable $\vartheta = (1 - \upsilon)\pi + \upsilon\theta_f - 2\tau$:
\begin{equation*}
    \tilde{\mathcal{B}}^\upsilon_\mathcal{P} = \left\{\bar{\boldsymbol{z}}_{\mathcal{B}_\mathcal{P}^\upsilon}(\tau, \vartheta):\: \vartheta \in (0, 2\pi),\: \tau \in \left(0, \pi - \frac\vartheta2\right)\right\},
\end{equation*}
where
\begin{multline*}
    \bar{\boldsymbol{z}}_{\mathcal{B}_\mathcal{P}^\upsilon}(\tau, \vartheta) \overset{\mathrm{def}}{=} \boldsymbol{z}_{\mathcal{B}_\mathcal{P}^\upsilon}(\tau, (1 - \upsilon)\pi + \upsilon\vartheta + 2\upsilon\tau)\\
    =
    \begin{bmatrix}
        &-\upsilon\left(\ell\sin\frac\vartheta2 + \cos\vartheta - \cos(\tau + \vartheta) + 1 - \cos\tau\right)\\
        &-\ell\cos\frac\vartheta2 + \sin\vartheta - \sin(\tau + \vartheta) + \sin\tau\\
        &(1 - \upsilon)\pi + \upsilon\vartheta
    \end{bmatrix}.
\end{multline*}
This parametrization is more convenient since the third component of the state vector depends only on $\vartheta$.

Let us now turn to the emanation from $\mathcal{BUP}_0$. The candidate optimal control laws $\tilde{u}(\cdot)$ and $\tilde{v}(\cdot)$ are equal to
\begin{equation}\label{eq:BUP0_controls}
    \begin{aligned}
        &\tilde{u}(\boldsymbol{z}_{\mathcal{BUP}_0}(\varphi_f), \boldsymbol{n}_{\mathcal{BUP}_0}(\varphi_f)) =
        \begin{cases}
            &\mathrm{sgn}\:\sin\varphi_f, \quad \varphi_f \neq 0_{\mathbb{S}};\\
            &u \in \{-1, 0, +1\}, \quad \varphi_f = 0_{\mathbb{S}},
        \end{cases}\\
        &\tilde{v}(\boldsymbol{z}_{\mathcal{BUP}_0}(\varphi_f), \boldsymbol{n}_{\mathcal{BUP}_0}(\varphi_f)) =
        \begin{cases}
            &\mathrm{sgn}\:\sin\varphi_f, \quad \varphi_f \neq 0_{\mathbb{S}};\\
            &v \in \{-1, 0, +1\}, \quad \varphi_f = 0_{\mathbb{S}}.
        \end{cases}
    \end{aligned}
\end{equation}
If $\varphi_f \neq 0_\mathbb{S}$, then the controls are equal to each other. Substituting $u = v = \upsilon$ and $\tilde{\boldsymbol{z}} = \boldsymbol{z}_{\mathcal{BUP}_0}(\varphi_f)$ into~\eqref{eq:radial_distance} yields $r_{\upsilon, \upsilon}(\tau; \boldsymbol{z}_{\mathcal{BUP}_0}(\varphi_f)) = \ell$. It means that the state vector remains on $\mathcal{BUP}_0$, and emanation doesn't proceed.

Next, if $\varphi_f = 0_\mathbb{S}$, then, according to~\eqref{eq:BUP0_controls}, we must analyze nine scenarios $u, v \in \{-1, 0, +1\}$. By Taylor's formula for the expression~\eqref{eq:radial_distance}, the distance between players is the following:
\begin{equation}\label{eq:radial_distance_taylor}
    r_{u, v}(\tau; \boldsymbol{z}_{\mathcal{BUP}_0}(0_\mathbb{S})) = \ell + \frac{v^2 - u^2}6\tau^3 + \frac{(u - v)^2}{8\ell}\tau^4 + o(\tau^4).
\end{equation}
If $u = v \in \{-1, 0, +1\}$, then, as in the previous case, $r_{v, v}(\tau; \boldsymbol{z}_{\mathcal{BUP}_0}(0_\mathbb{S})) = \ell$.

Next, if $u \in \mathbb{B}$ and $v = 0$, then, according to~\eqref{eq:radial_distance_taylor}, the distance becomes less than the capture radius at the start of emanation. The region inside the capture circle doesn't belong to the game set and, in consequence, this case can be missed.

Now, let $v = -u = \upsilon \in \mathbb{B}$. We associate the notions $\mathcal{B}_\mathcal{PL}^\upsilon$ and $\tilde{\mathcal{B}}_\mathcal{PL}^\upsilon$ with the part of barrier $\mathcal{B}$, that emanates from $\begin{bmatrix} 0 & \ell & 0_\mathbb{S} \end{bmatrix}^\mathrm{T} \in \mathcal{BUP}_0$ by analogy with $\mathcal{B}_\mathcal{P}^\upsilon$ and $\tilde{\mathcal{B}}_\mathcal{P}^\upsilon$ notions. Note that $\tilde{\mathcal{B}}_\mathcal{PL}^\upsilon$ is a line while $\tilde{\mathcal{B}}_\mathcal{P}^\upsilon$ is a surface. In the same manner, we can describe the points on $\tilde{\mathcal{B}}_\mathcal{PL}^\upsilon$
\begin{equation*}
    \boldsymbol{z}_{\mathcal{B}_\mathcal{PL}^\upsilon}(\tau) \overset{\mathrm{def}}{=} \boldsymbol{z}_{-\upsilon, \upsilon}(\tau; \boldsymbol{z}_{\mathcal{BUP}_0}(0_\mathbb{S})) =
    \begin{bmatrix}
        -\upsilon\left(\ell\sin\tau - 2\cos\tau + 1 + \cos2\tau\right)\\
        \ell\cos\tau + 2\sin{\tau} - \sin2\tau\\
        (1 + \upsilon)\pi - 2\upsilon\tau
    \end{bmatrix}
\end{equation*}
and the corresponding normal vector
\begin{equation*}
    \nu_{\boldsymbol{z}, \mathcal{B}_\mathcal{PL}^\upsilon}(\tau) \overset{\mathrm{def}}{=} \nu^{-\upsilon, \upsilon}_{\boldsymbol{z}}(\tau; \boldsymbol{z}_{\mathcal{BUP}_0}(0_\mathbb{S}), \boldsymbol{n}_{\mathcal{BUP}_0}(0_\mathbb{S})).
\end{equation*}
The analysis of the switch functions 
\begin{align*}
    s_P(\boldsymbol{z}_{\mathcal{B}_\mathcal{PL}^\upsilon}(\tau), \nu_{\boldsymbol{z}, \mathcal{B}_\mathcal{PL}^\upsilon}(\tau)) = \upsilon(\cos\tau - 1) = 0,\\
    s_E(\boldsymbol{z}_{\mathcal{B}_\mathcal{PL}^\upsilon}(\tau), \nu_{\boldsymbol{z}, \mathcal{B}_\mathcal{PL}^\upsilon}(\tau)) = \upsilon(1 - \cos\tau) = 0
\end{align*}
gives the minimal switch time $\tau = 2\pi$. The third component of $\boldsymbol{z}_{\mathcal{B}_\mathcal{PL}^\upsilon}(\tau)$ vanishes at $\tau = \pi$. Therefore, it may be concluded that 
\begin{equation}\label{eq:BPLv_definition}
    \tilde{\mathcal{B}}_\mathcal{PL}^\upsilon \overset{\mathrm{def}}{=} \{\boldsymbol{z}_{\mathcal{B}_\mathcal{PL}^\upsilon}(\tau):\: \tau \in (0, \pi)\}. 
\end{equation}
Applying a new parametrization $\tau = \pi - \vartheta/2$ we can rewrite~\eqref{eq:BPLv_definition} as
\begin{equation*}
    \tilde{\mathcal{B}}_\mathcal{PL}^\upsilon = \{\bar{\boldsymbol{z}}_{\mathcal{B}_\mathcal{PL}^\upsilon}(\vartheta):\: \vartheta \in (0, 2\pi)\}, 
\end{equation*}
where
\begin{equation*}
    \bar{\boldsymbol{z}}_{\mathcal{B}_\mathcal{PL}^\upsilon}(\vartheta) \overset{\mathrm{def}}{=} \boldsymbol{z}_{\mathcal{B}_\mathcal{PL}^\upsilon}\left(\pi - \frac\vartheta2\right) =
    \begin{bmatrix}
        -\upsilon\left(\ell\sin\frac\vartheta2 + 2\cos\frac\vartheta2 + 1 + \cos\vartheta\right)\\
        -\ell\cos\frac\vartheta2 + 2\sin\frac\vartheta2 + \sin\vartheta\\
        (1 - \upsilon)\pi + \upsilon\vartheta
    \end{bmatrix}.
\end{equation*}
An easy computation shows that $\bar{\boldsymbol{z}}_{\mathcal{B}_\mathcal{P}^\upsilon}(\vartheta, \pi - \vartheta/2) = \bar{\boldsymbol{z}}_{\mathcal{B}_\mathcal{PL}^\upsilon}(\vartheta)$ for all $\vartheta \in (0, 2\pi)$, which implies that $\tilde{\mathcal{B}}_\mathcal{PL}^\upsilon$ is a part of the boundary of $\tilde{\mathcal{B}}_\mathcal{P}^\upsilon$.

Finally, consider the case $u = 0$, $v = \upsilon \in \mathbb{B}$. We will see later that this case corresponds to emanation of a universal line on the barrier. The notions $\mathcal{B}_\mathcal{UL}^\upsilon$, $\tilde{\mathcal{B}}_\mathcal{UL}^\upsilon$ are associated with the emanation. The points of $\tilde{\mathcal{B}}_\mathcal{UL}^\upsilon$ are given by
\begin{equation*}
    \boldsymbol{z}_{\mathcal{B}_\mathcal{UL}^\upsilon}(\tau) \overset{\mathrm{def}}{=} \boldsymbol{z}_{0, \upsilon}(\tau; \boldsymbol{z}_{\mathcal{BUP}_0}(0_\mathbb{S})) = 
    \begin{bmatrix}
        \upsilon(1 - \cos\tau)\\
        \ell + \tau - \sin{\tau}\\
        (1 + \upsilon)\pi - \upsilon\tau
    \end{bmatrix}.
\end{equation*}
The normal vector is defined by
\begin{equation}\label{eq:BUL_normal}
    \nu_{\boldsymbol{z}, \mathcal{B}_\mathcal{UL}^\upsilon}(\tau) \overset{\mathrm{def}}{=} \nu^{0, \upsilon}_{\boldsymbol{z}}(\tau; \boldsymbol{z}_{\mathcal{BUP}_0}(0_\mathbb{S}), \boldsymbol{n}_{\mathcal{BUP}_0}(0_\mathbb{S})) =
    \begin{bmatrix}
        0\\
        1\\
        \upsilon(1 - \cos\tau)
    \end{bmatrix}.
\end{equation}
Let us compute the switch functions on $\tilde{\mathcal{B}}_\mathcal{UL}^\upsilon$:
\begin{equation*}
    \begin{aligned}
        &s_P(\boldsymbol{z}_{\mathcal{B}_\mathcal{UL}^\upsilon}(\tau), \nu_{\boldsymbol{z}, \mathcal{B}_\mathcal{UL}^\upsilon}(\tau)) = 0,\\
        &s_E(\boldsymbol{z}_{\mathcal{B}_\mathcal{UL}^\upsilon}(\tau), \nu_{\boldsymbol{z}, \mathcal{B}_\mathcal{UL}^\upsilon}(\tau)) = \upsilon(1 - \cos\tau).
    \end{aligned}
\end{equation*}
Since the first component of~\eqref{eq:BUL_normal} is equal to zero and the second component is non-negative, the emanation with $u \in \{-1, 0, +1\}$ is not prohibited on $\tilde{\mathcal{B}}_\mathcal{UL}^\upsilon$. The minimal switch time of the evader switch function is $\tau = 2\pi$. Thus
\begin{equation*}
    \tilde{\mathcal{B}}_\mathcal{UL}^\upsilon \overset{\mathrm{def}}{=} \{\boldsymbol{z}_{\mathcal{B}_\mathcal{UL}^\upsilon}(\tau):\: \tau \in (0, 2\pi)\}.
\end{equation*}
For uniformity of notation, we define $\bar{\boldsymbol{z}}_{\mathcal{B}_\mathcal{UL}^\upsilon}(\vartheta) = \boldsymbol{z}_{\mathcal{B}_\mathcal{UL}^\upsilon}\left(\vartheta\right)$.

Thus, the primary emanation of trajectories from the $\mathcal{BUP}$ is exhausted by two surfaces -- $\tilde{\mathcal{B}}_\mathcal{P}^\upsilon$ and four lines -- $\tilde{\mathcal{B}}_\mathcal{PL}^\upsilon$, $\tilde{\mathcal{B}}_\mathcal{UL}^\upsilon$ ($\upsilon \in \mathbb{B}$). At the same time, each line $\tilde{\mathcal{B}}_\mathcal{PL}^\upsilon$ is the boundary of surface $\tilde{\mathcal{B}}_\mathcal{P}^\upsilon$. Next, we describe the tributaries of the $\tilde{\mathcal{B}}_\mathcal{UL}^\upsilon$, which consist of retrograde trajectories emanating in a retrograde sense from this line.

\subsection{Tributaries emanation}

Since the tributaries emanate from $\tilde{\mathcal{B}}_\mathcal{UL}^\upsilon$ with two possible values of the pursuer's control $u \in \mathbb{B}$, we explore the emanation in two scenarios. In the first scenario, we suppose that players' controls are equal to each other and the players turn in the same direction. In the second case, the players turn in the opposite directions, i.e., the controls of the players are different.

We associate notions $\mathcal{B}_\mathcal{TS}^\upsilon$ and $\tilde{\mathcal{B}}_\mathcal{TS}^\upsilon$ with the emanation of tributaries with the same values of players controls. Substituting $u = v = \upsilon$ into the solution of equations of motion for fixed values of controls we have
\begin{equation*}
    \boldsymbol{z}_{\mathcal{B}_\mathcal{TS}^\upsilon}(\tau_1, \tau_2) \overset{\mathrm{def}}{=} \boldsymbol{z}_{\upsilon, \upsilon}(\tau_2; \boldsymbol{z}_{\mathcal{B}_\mathcal{UL}^\upsilon}(\tau_1)) =
    \begin{bmatrix}
        \upsilon\left((\ell + \tau_1)\sin\tau_2 + 1 - \cos\tau_1\right)\\
        (\ell + \tau_1)\cos\tau_2 - \sin\tau_1\\
        (1 + \upsilon)\pi - \upsilon\tau_1
    \end{bmatrix}.
\end{equation*}
Here, $\tau_1 \in (0, 2\pi)$ is the time spent on $\tilde{\mathcal{B}}_\mathcal{UL}^\upsilon$, and $\tau_2$ is spent on the tributary path. The normal vector is defined by
\begin{equation*}
    \nu_{\boldsymbol{z}, \mathcal{B}_\mathcal{TS}^\upsilon}(\tau_1, \tau_2) \overset{\mathrm{def}}{=} \nu_{\boldsymbol{z}}^{\upsilon, \upsilon}(\tau_2; \boldsymbol{z}_{\mathcal{B}_\mathcal{UL}^\upsilon}(\tau_1), \nu_{\boldsymbol{z}, \mathcal{B}_\mathcal{UL}^\upsilon}(\tau_1)).
\end{equation*}
Computing the switch functions on the $\tilde{\mathcal{B}}_\mathcal{TS}^\upsilon$ gives
\begin{equation*}
    \begin{aligned}
        &s_P(\boldsymbol{z}_{\mathcal{B}_\mathcal{TS}^\upsilon}(\tau_1, \tau_2), \nu_{\boldsymbol{z}, \mathcal{B}_\mathcal{TS}^\upsilon}(\tau_1, \tau_2)) = \upsilon(1 - \cos\tau_2) = 0,\\
        &s_E(\boldsymbol{z}_{\mathcal{B}_\mathcal{TS}^\upsilon}(\tau_1, \tau_2), \nu_{\boldsymbol{z}, \mathcal{B}_\mathcal{TS}^\upsilon}(\tau_1, \tau_2)) = \upsilon\left(1 - \cos\left(\tau_1 + \tau_2\right)\right) = 0.
    \end{aligned}
\end{equation*}
The minimal positive solution of these equations $\tau_2 = 2\pi - \tau_1$ corresponds to zero of the evader's switch function. Hence
\begin{equation*}
        \tilde{\mathcal{B}}_\mathcal{TS}^\upsilon \overset{\mathrm{def}}{=} \left\{\boldsymbol{z}_{\mathcal{B}_\mathcal{TS}^\upsilon}\left(\tau_1, \tau_2\right):\: \tau_1 \in \left(0, 2\pi\right), \: \tau_2 \in \left(0, 2\pi - \tau_1\right)\right\}.
\end{equation*}
Using a new parametrization $\tau_1 = \vartheta$, $\tau_2 = \tau - \vartheta$ we can rewrite the above as
\begin{equation*}
    \tilde{\mathcal{B}}_\mathcal{TS}^\upsilon = \left\{\bar{\boldsymbol{z}}_{\mathcal{B}_\mathcal{TS}^\upsilon}\left(\tau, \vartheta\right):\: \vartheta \in \left(0, 2\pi\right), \:\tau \in \left(\vartheta, 2\pi\right)\right\},
\end{equation*}
where
\begin{equation*}
    \bar{\boldsymbol{z}}_{\mathcal{B}_\mathcal{TS}^\upsilon}(\tau, \vartheta) \overset{\mathrm{def}}{=} \boldsymbol{z}_{\mathcal{B}_\mathcal{TS}^\upsilon}\left(\vartheta, \tau - \vartheta\right) =
    \begin{bmatrix}
        \upsilon\left((\ell + \vartheta)\sin(\tau - \vartheta) + 1 - \cos\vartheta\right)\\
        (\ell + \vartheta)\cos(\tau - \vartheta) - \sin\vartheta\\
        (1 + \upsilon)\pi - \upsilon\vartheta
    \end{bmatrix}.
\end{equation*}
Note that $\tau = \tau_1 + \tau_2$ has a sense of the time-to-go to the $\mathcal{BUP}_0$.

Next, we associate notions $\mathcal{B}_\mathcal{TD}^\upsilon$ and $\tilde{\mathcal{B}}_\mathcal{TD}^\upsilon$ with the emanation of tributaries with the different values of players' controls. Substituting $v = -u = \upsilon$ into the solution of equations of motion we obtain
\begin{multline*}
    \boldsymbol{z}_{\mathcal{B}_\mathcal{TD}^\upsilon}(\tau_1, \tau_2) \overset{\mathrm{def}}{=} \boldsymbol{z}_{-\upsilon, \upsilon}(\tau_2; \boldsymbol{z}_{\mathcal{B}_\mathcal{UL}^\upsilon}(\tau_1))\\
    =
    \begin{bmatrix}
        \upsilon\left(-(\ell + \tau_1)\sin\tau_2 + 2\cos\tau_2 - \cos(\tau_1 + 2\tau_2) - 1\right)\\
        (\ell + \tau_1)\cos\tau_2 + 2\sin\tau_2 - \sin(\tau_1 + 2\tau_2)\\
        (1 + \upsilon)\pi - \upsilon(\tau_1 + 2\tau_2)
    \end{bmatrix}.
\end{multline*}
The variables $\tau_1$, $\tau_2$ have the same sense as in the previous case. The normal vector is
\begin{equation*}
    \nu_{\boldsymbol{z}, \mathcal{B}_\mathcal{TD}^\upsilon}(\tau_1, \tau_2) \overset{\mathrm{def}}{=} \nu_{\boldsymbol{z}}^{-\upsilon, \upsilon}(\tau_2; \boldsymbol{z}_{\mathcal{B}_\mathcal{UL}^\upsilon}(\tau_1), \nu_{\boldsymbol{z}, \mathcal{B}_\mathcal{UL}^\upsilon}(\tau_1)).
\end{equation*}
Calculating the switch functions gives
\begin{equation*}
    \begin{aligned}
        &s_P(\boldsymbol{z}_{\mathcal{B}_\mathcal{TD}^\upsilon}(\tau_1, \tau_2), \nu_{\boldsymbol{z}, \mathcal{B}_\mathcal{TD}^\upsilon}(\tau_1, \tau_2)) = \upsilon(\cos\tau_2 - 1) = 0,\\
        &s_E(\boldsymbol{z}_{\mathcal{B}_\mathcal{TD}^\upsilon}(\tau_1, \tau_2), \nu_{\boldsymbol{z}, \mathcal{B}_\mathcal{TD}^\upsilon}(\tau_1, \tau_2)) = \upsilon\left(1 - \cos\left(\tau_1 + \tau_2\right)\right) = 0.
    \end{aligned}
\end{equation*}
The minimal positive $\tau_2 = 2\pi - \tau_1$ also corresponds to zero of the evader's switch function, but the third component of $\boldsymbol{z}_{\mathcal{B}_\mathcal{TD}^\upsilon}(\tau_1, \tau_2)$ vanishes earlier at $\tau_2 = \pi - \tau_1/2$. Therefore, we set
\begin{equation*}
    \tilde{\mathcal{B}}_\mathcal{TD}^\upsilon \overset{\mathrm{def}}{=} \left\{\boldsymbol{z}_{\mathcal{B}_\mathcal{TD}^\upsilon}\left(\tau_1, \tau_2\right):\: \tau_1 \in \left(0, 2\pi\right), \: \tau_2 \in \left(0, \pi - \frac{\tau_1}2\right)\right\}.
\end{equation*}
Using a new parametrization $\tau_1 = 2\tau - \vartheta$, $\tau_2 = \vartheta - \tau$ we obtain 
\begin{equation*}
    \tilde{\mathcal{B}}_\mathcal{TD}^\upsilon = \left\{\bar{\boldsymbol{z}}_{\mathcal{B}_\mathcal{TD}^\upsilon}\left(\tau, \vartheta\right):\: \vartheta \in \left(0, 2\pi\right), \:\tau \in \left(\frac\vartheta2, \vartheta\right)\right\},
\end{equation*}
where
\begin{multline*}
    \bar{\boldsymbol{z}}_{\mathcal{B}_\mathcal{TD}^\upsilon}(\tau, \vartheta) \overset{\mathrm{def}}{=} \boldsymbol{z}_{\mathcal{B}_\mathcal{TD}^\upsilon}\left(2\tau - \vartheta, \vartheta - \tau\right)\\
    =
    \begin{bmatrix}
        \upsilon\left((\ell + 2\tau - \vartheta)\sin(\tau - \vartheta) + 2\cos(\tau - \vartheta) - 1 - \cos\vartheta\right)\\
        (\ell + 2\tau - \vartheta)\cos(\tau - \vartheta) - 2\sin(\tau - \vartheta) - \sin\vartheta\\
        (1 + \upsilon)\pi - \upsilon\vartheta
    \end{bmatrix}.
\end{multline*}
An easy computation shows that $\bar{\boldsymbol{z}}_{\mathcal{B}_\mathcal{TD}^\upsilon}(\vartheta/2, \vartheta) = \bar{\boldsymbol{z}}_{\mathcal{B}_\mathcal{PL}^\upsilon}(2\pi-\vartheta)$ for all $\vartheta \in (0, 2\pi)$, which proves that line $\tilde{\mathcal{B}}_\mathcal{PL}^\upsilon$ is a part of the boundary of the surface $\tilde{\mathcal{B}}_\mathcal{TD}^{\upsilon}$.

\subsection{Barrier self-intersections}

We now turn to the problem of determining "redundant" parts of obtained semipermeable surfaces. These "redundant" parts are semipermeable but they are not involved in barrier $\mathcal{B}$. The GTIC has a symmetry $(x, y, \theta, u, v) \leftrightarrow (-x, y, 2\pi - \theta, -u, -v)$, i.e. the equations of motion~\eqref{eq:dynamics} and the capture set $\mathcal{C}$ hold when we use this replacement of variables. This fact allows exploring the GTIC only for $\theta \in (0, \pi]$. For the rest of the values $\theta \in (\pi, 2\pi)$, we can use this symmetry.

\begin{figure}
    \begin{center}
        \includegraphics[width=0.434\linewidth]{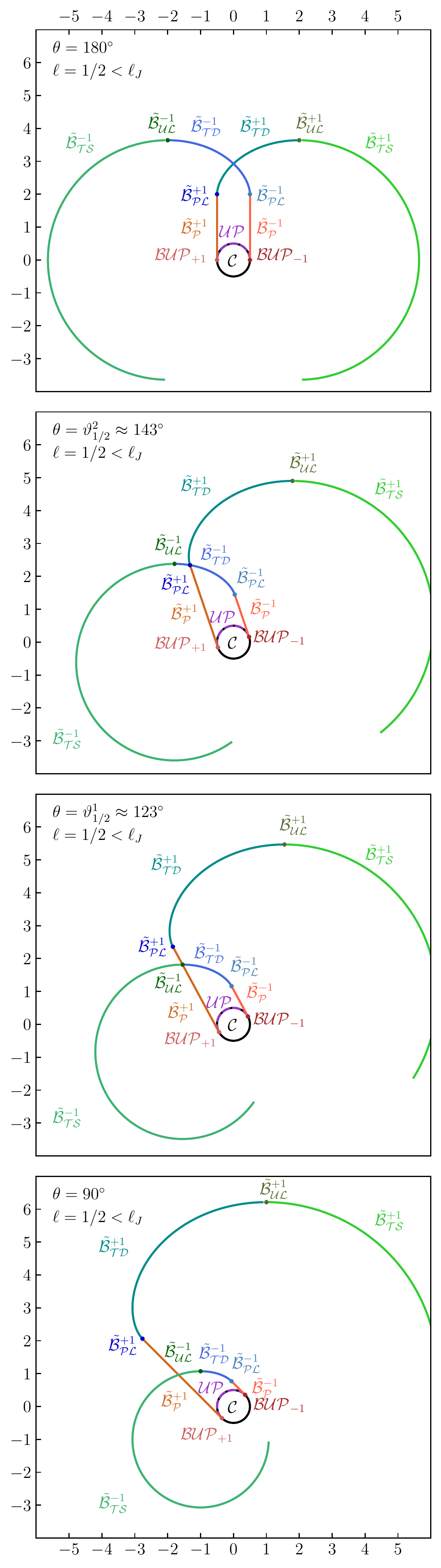}
        \includegraphics[width=0.46\linewidth]{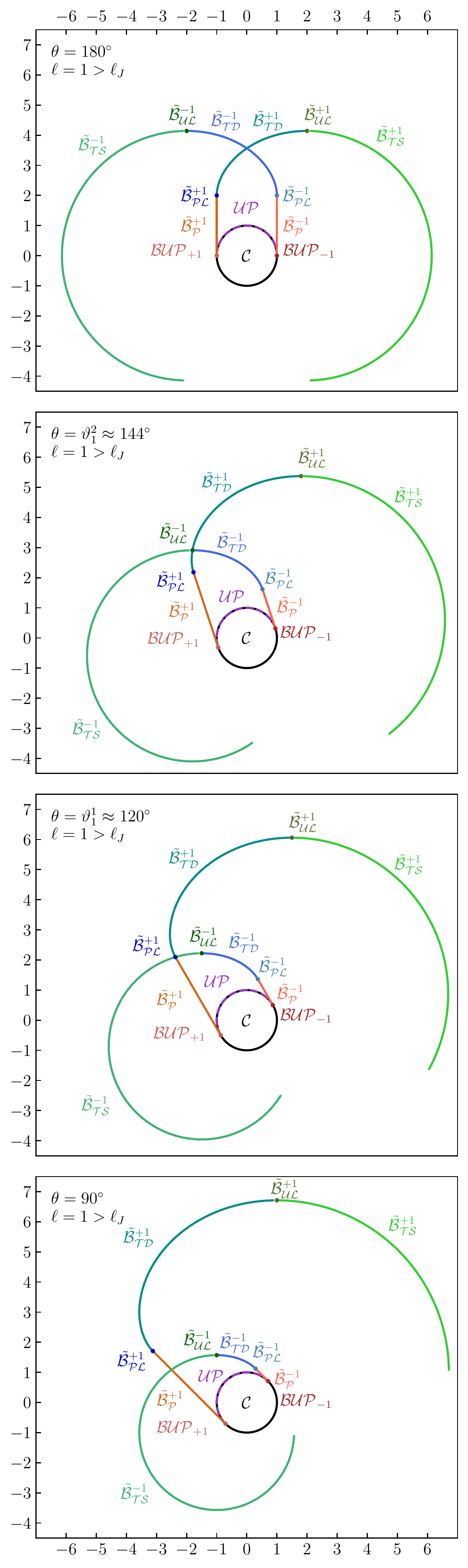}
        \caption{Cross sections of the emanated surfaces for the "small" and "large" capture radii}
        \label{fig:small_and_large_radius_slices}
    \end{center}
\end{figure}

\begin{figure}
    \begin{center}
        \includegraphics[width=\linewidth]{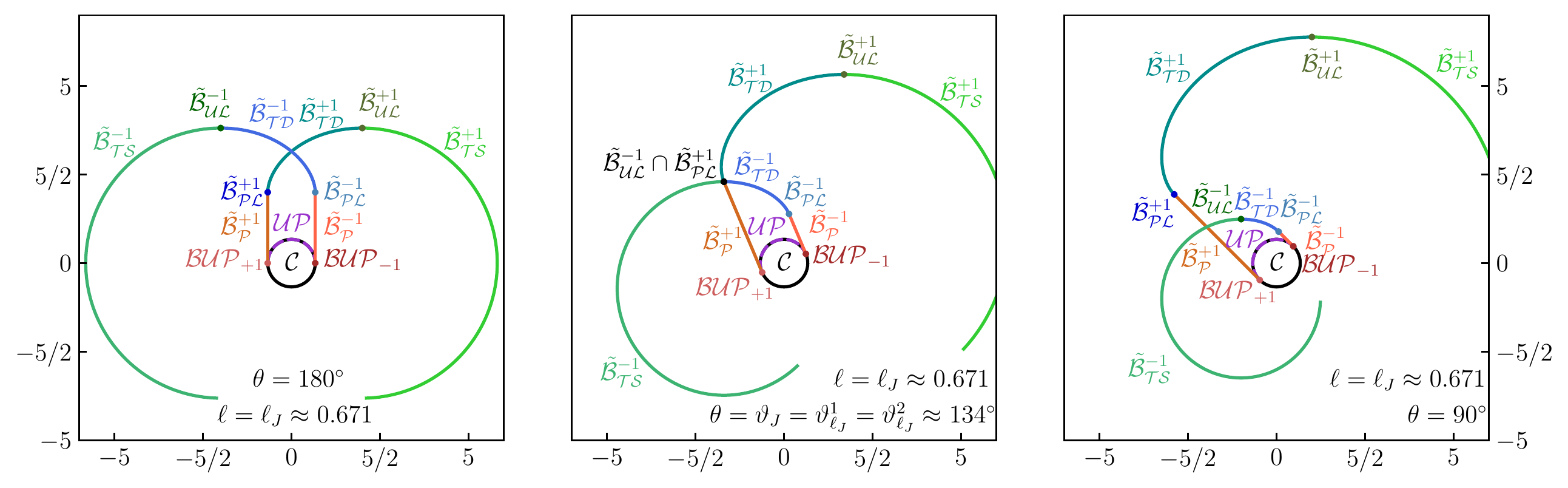}
        \caption{Cross sections of the emanated surfaces for the "medium" capture radius $\ell = \ell_J$}
        \label{fig:middle_radius_slices}
    \end{center}
\end{figure}

Graphical analysis of the $\theta$-slices of emanating parts (see Fig.~\ref{fig:small_and_large_radius_slices}) shows that for small enough values of the capture radius $\ell$ (e.g., for $\ell = 1/2$), the surface $\tilde{\mathcal{B}}_\mathcal{P}^\upsilon$ crosses the $\tilde{\mathcal{B}}_\mathcal{TD}^{-\upsilon}$, but the surface $\tilde{\mathcal{B}}_\mathcal{TD}^\upsilon$ doesn't cross the $\tilde{\mathcal{B}}_\mathcal{TS}^{-\upsilon}$. At the same time, however, for large enough values of $\ell$ (e.g., for $\ell = 1$), we observe an opposite behavior. This fact suggests that there exists an intermediate value of $\ell = \ell_J$ when the surface $\tilde{\mathcal{B}}_\mathcal{P}^\upsilon$ doesn't cross the $\tilde{\mathcal{B}}_\mathcal{TD}^{-\upsilon}$ and the $\tilde{\mathcal{B}}_\mathcal{TD}^\upsilon$ doesn't cross the $\tilde{\mathcal{B}}_\mathcal{TS}^{-\upsilon}$, but the boundaries of these surfaces share the common point for some $\theta = (1 - \upsilon)\pi + \upsilon\vartheta_J$, where $\vartheta_J \in (0, \pi)$ (see Fig.~\ref{fig:middle_radius_slices}).


\begin{lemma}
    There is only one value of the capture radius $\ell = \ell_J$ when line $\tilde{\mathcal{B}}_\mathcal{PL}^\upsilon$ crosses line $\tilde{\mathcal{B}}_\mathcal{UL}^{-\upsilon}$.
\end{lemma}
For future purposes, we introduce a special rotation matrix
\begin{equation*}
    \boldsymbol{R}_\upsilon(\alpha) \overset{\mathrm{def}}{=}
    \begin{bmatrix}
        -\upsilon\cos\alpha & \sin\alpha & 0 \\
        \upsilon\sin\alpha & \cos\alpha & 0 \\
        0 & 0 & 1
    \end{bmatrix}, \quad \alpha \in \mathbb{S}, \quad \upsilon \in \mathbb{B}.
\end{equation*}
Note that for all $\upsilon \in \mathbb{B}$ and $\alpha \in \mathbb{S}$ the matrix is non-singular. Also, we will use the following functions:
\begin{equation*}
    \xi_\ell(\alpha) \overset{\mathrm{def}}{=} (\ell + \alpha)\sin\frac\alpha2 + 2\cos\frac\alpha2, \quad \eta_\ell(\alpha) \overset{\mathrm{def}}{=} (\ell + \alpha)\cos\frac\alpha2 - 2\sin\frac\alpha2.
\end{equation*}

\begin{proof}
    Let us first prove that the system
    \begin{equation*}
        \left\{
        \begin{aligned}
            &\bar{\boldsymbol{z}}_{\mathcal{B}_\mathcal{UL}^{-\upsilon}}(\vartheta) - \bar{\boldsymbol{z}}_{\mathcal{B}_\mathcal{PL}^\upsilon}(\vartheta) = \boldsymbol{0},\\
            &\vartheta \in (0, 2\pi), \quad \ell \in \mathbb{R}^+
        \end{aligned}
        \right.
    \end{equation*}
    has only one solution for each $\upsilon \in \mathbb{B}$. Applying $\boldsymbol{R}_\upsilon(\vartheta/2)$ we can rewrite the above as
    \begin{equation*}
        \left\{
        \begin{aligned}
            &\boldsymbol{R}_\upsilon\left(\frac\vartheta2\right)\left(\bar{\boldsymbol{z}}_{\mathcal{B}_\mathcal{UL}^{-\upsilon}}(\vartheta) - \bar{\boldsymbol{z}}_{\mathcal{B}_\mathcal{PL}^\upsilon}(\vartheta)\right) =
            \begin{bmatrix}
                \xi_\ell(\vartheta) - 2 - 4\cos\frac\vartheta2\\
                \eta_\ell(\vartheta) + \ell\\
                0_\mathbb{S}
            \end{bmatrix} = \boldsymbol{0},\\
            &\vartheta \in (0, 2\pi), \quad \ell \in \mathbb{R}^+\\
        \end{aligned}
        \right.
    \end{equation*}
    Eliminating the variable $\ell$ from the system yields
    \begin{equation}\label{eq:thetaJ}
        \vartheta - 4\left(1 + \cos\frac\vartheta2\right)\cot\frac\vartheta2 = 0, \quad \vartheta \in (0, 2\pi).
    \end{equation}
    The left part of \eqref{eq:thetaJ} approaches to $-\infty$ as $\vartheta \to +0$ and $2\pi$ as $\vartheta \to 2\pi - 0$. The derivative of the left part equals
    \begin{equation*}
         1 + 2\left(\cos\frac\vartheta2 + \frac{1 + \cos\frac\vartheta2}{\sin^2\frac\vartheta2}\right) = \frac{\left(1 + \cos\frac\vartheta2\right)^2\left(3 - 2\cos\frac\vartheta2\right)}{\sin^2\frac\vartheta2} > 0, \quad \vartheta \in (0, 2\pi).
    \end{equation*}
    Thus, the left part of \eqref{eq:thetaJ} is a continuous monotonic function. It has opposite signs for $\vartheta \to +0$ and $\vartheta \to 2\pi - 0$. Hence, the equation \eqref{eq:thetaJ} has only one root. \qed
\end{proof}

Applying Newton's method\footnote{An initial guess is $\vartheta_0 = 2$, for example.} gives the value $\vartheta_J \approx 2.343 \approx 134^\circ$ for the root. The corresponding capture radius is given by
\begin{equation*}
    \ell_J \overset{\mathrm{def}}{=} -2\frac{\cos\frac{\vartheta_J}2 + \cos\vartheta_J}{\sin\frac{\vartheta_J}2} \approx 0.671.
\end{equation*}

The intersection of the semipermeable surfaces constituting the barrier indicates the presence of a dispersal line on the barrier. Analysis of $\theta$-slices in Figs.~\ref{fig:small_and_large_radius_slices}-\ref{fig:middle_radius_slices} shows that the dispersal line is built in different ways, depending on the value of the capture radius $\ell$. At the same time, for values of $\theta$ close to $0_\mathbb{S}$ for all values of the capture radius, there is an intersection $\tilde{\mathcal{B}}_\mathcal{P}^\upsilon \cap \tilde{\mathcal{B}}_\mathcal {TS}^{-\upsilon}$, and for $\theta$ close to $\pi$ there is an intersection $\tilde{\mathcal{B}}_\mathcal{TD}^\upsilon \cap \tilde{\mathcal{B}}_ \mathcal{TD}^{-\upsilon}$. We will use notation $\mathcal{B}_\mathcal{DL}$ for the dispersal line. Taking into account all of these facts, we distinguish three cases of constructing a dispersal line:
\begin{itemize}
    \item the "small" capture radius $\ell \in (0, \ell_J)$
    \begin{multline*}
        \mathcal{B}_\mathcal{DL} \overset{\mathrm{def}}{=} \bigcup\limits_{\upsilon \in \mathbb{B}} \left[\left(\tilde{\mathcal{B}}_\mathcal{P}^\upsilon \cap \tilde{\mathcal{B}}_\mathcal{TS}^{-\upsilon}\right) \cup \left( \tilde{\mathcal{B}}_\mathcal{P}^\upsilon \cap \tilde{\mathcal{B}}_\mathcal{UL}^{-\upsilon}\right) \cup \left(\tilde{\mathcal{B}}_\mathcal{P}^\upsilon \cap \tilde{\mathcal{B}}_\mathcal{TD}^{-\upsilon}\right) \right.\\
        \left. \cup \left(\tilde{\mathcal{B}}_\mathcal{PL}^\upsilon \cap \tilde{\mathcal{B}}_\mathcal{TD}^{-\upsilon}\right) \cup \left(\tilde{\mathcal{B}}_\mathcal{TD}^\upsilon \cap \tilde{\mathcal{B}}_\mathcal{TD}^{-\upsilon}\right)\right];
    \end{multline*}
    \item the "medium" capture radius $\ell = \ell_J$
    \begin{equation*}
        \mathcal{B}_\mathcal{DL} \overset{\mathrm{def}}{=} \bigcup\limits_{\upsilon \in \mathbb{B}} \left[\left(\tilde{\mathcal{B}}_\mathcal{P}^\upsilon \cap \tilde{\mathcal{B}}_\mathcal{TS}^{-\upsilon}\right) \cup \left(\tilde{\mathcal{B}}_\mathcal{PL}^\upsilon \cap \tilde{\mathcal{B}}_\mathcal{UL}^{-\upsilon}\right) \cup \left(\tilde{\mathcal{B}}_\mathcal{TD}^\upsilon \cap \tilde{\mathcal{B}}_\mathcal{TD}^{-\upsilon}\right)\right];
    \end{equation*}
    \item the "large" capture radius $\ell \in (\ell_J, +\infty)$
    \begin{multline*}
            \mathcal{B}_\mathcal{DL} \overset{\mathrm{def}}{=} \bigcup\limits_{\upsilon \in \mathbb{B}}\left[\left(\tilde{\mathcal{B}}_\mathcal{P}^\upsilon \cap \tilde{\mathcal{B}}_\mathcal{TS}^{-\upsilon}\right) \cup \left( \tilde{\mathcal{B}}_\mathcal{PL}^\upsilon \cap \tilde{\mathcal{B}}_\mathcal{TS}^{-\upsilon}\right) \cup \left(\tilde{\mathcal{B}}_\mathcal{TD}^\upsilon \cap \tilde{\mathcal{B}}_\mathcal{TS}^{-\upsilon}\right) \right.\\
            \left. \cup \left(\tilde{\mathcal{B}}_\mathcal{TD}^\upsilon \cap \tilde{\mathcal{B}}_\mathcal{UL}^{-\upsilon}\right) \cup \left(\tilde{\mathcal{B}}_\mathcal{TD}^\upsilon \cap \tilde{\mathcal{B}}_\mathcal{TD}^{-\upsilon}\right)\right].
    \end{multline*}
\end{itemize}

To determine which range of $\theta$ corresponds to the intersection $\tilde{\mathcal{B}}_\mathcal{P}^\upsilon \cap \tilde{\mathcal{B}}_\mathcal{TS}^{ -\upsilon}$, we must find the critical value of the angle at which the universal line $\tilde{\mathcal{B}}_\mathcal{UL}^{-\upsilon}$ crosses the surface $\tilde{\mathcal {B}}_\mathcal{P}^\upsilon$ if $\ell$ is "small". If $\ell$ is "large", then the critical value of the angle corresponds to the intersection of line $\tilde{\mathcal{B}}_ \mathcal{PL}^\upsilon$ with the surface $\tilde{\mathcal{B}}_\mathcal{TS}^{-\upsilon}$. For the "medium" value of $\ell$, the critical value of the angle is $(1 - \upsilon)\pi + \upsilon\vartheta_J$. We denote by $\vartheta^1_\ell$ a dependence of the critical value of angle from the capture radius $\ell$ for $\upsilon = +1$. If $\upsilon = -1$, this angle equals $2\pi - \vartheta^1_\ell$ due to the symmetry of the problem.

The intersection of $\tilde{\mathcal{B}}_\mathcal{UL}^{-\upsilon}$ with $\tilde{\mathcal{B}}_\mathcal{P}^\upsilon$ is given by the system and interval constraints
\begin{equation*}
    \left\{
    \begin{aligned}
        &\bar{\boldsymbol{z}}_{\mathcal{B}_\mathcal{UL}^{-\upsilon}}(\vartheta) - \bar{\boldsymbol{z}}_{\mathcal{B}_\mathcal{P}^\upsilon}(\tau, \vartheta) = \boldsymbol{0},\\
        &\vartheta \in (0, 2\pi), \quad \tau \in \left(0, \pi - \frac\vartheta2\right), \quad \ell \in (0, \ell_J).
    \end{aligned}
    \right.
\end{equation*}
Transformation $\boldsymbol{R}_\upsilon(\vartheta/2)$ of the system yields 
\begin{equation*}
    \left\{
    \begin{aligned}
        &\xi_\ell(\vartheta) + 2\cos\left(\tau + \frac\vartheta2\right) - 4\cos\frac\vartheta2 = 0,\\
        &\eta_\ell(\vartheta) + \ell = 0,\\
        &\vartheta \in (0, 2\pi), \quad \tau \in \left(0, \pi - \frac\vartheta2\right), \quad \ell \in (0, \ell_J).
    \end{aligned}
    \right.
\end{equation*}
A unique solution of the system for the given ranges of $\vartheta$ and $\tau$ is
\begin{align*}
    &\tau = -\frac{w_\ell}2 + \arccos\left(2\cos\frac{w_\ell}2 - \frac{\sqrt{(\ell + w_\ell)^2 - \ell^2 + 4}}2\right),\\
    &\vartheta = w_\ell,
\end{align*}
where $w_\ell$ is a root of the transcendental equation\footnote{The root can be found by Newton's method with an initial guess $w_\ell^0 = \vartheta_J$, for example.}
\begin{equation}\label{eq:wl}
    \eta_\ell(w) + \ell = 0, \quad w \in (0, 2\pi), \quad \ell \in (0, \ell_J).
\end{equation}

The intersection of $\tilde{\mathcal{B}}_\mathcal{PL}^\upsilon$ with $\tilde{\mathcal{B}}_\mathcal{TS}^{-\upsilon}$ is determined by the system and interval constraints
\begin{equation*}
    \left\{
    \begin{aligned}
        &\bar{\boldsymbol{z}}_{\mathcal{B}_\mathcal{TS}^{-\upsilon}}(\tau, \vartheta) - \bar{\boldsymbol{z}}_{\mathcal{B}_\mathcal{PL}^\upsilon}(\vartheta)= \boldsymbol{0},\\
        &\vartheta \in (0, 2\pi), \quad \tau \in (\vartheta, 2\pi), \quad \ell \in (\ell_J, +\infty).
    \end{aligned}
    \right.
\end{equation*}
Transformation $\boldsymbol{R}_\upsilon(\vartheta/2)$ of the system yields
\begin{equation*}
    \left\{
    \begin{aligned}
        &(\ell + \vartheta)\sin\left(\tau - \frac\vartheta2\right) - 2\cos\frac\vartheta2 - 2 = 0,\\
        &(\ell + \vartheta)\cos\left(\tau - \frac\vartheta2\right) - 2\sin\frac\vartheta2 + \ell = 0,\\
        &\vartheta \in (0, 2\pi), \quad \tau \in (\vartheta, 2\pi), \quad \ell \in (\ell_J, +\infty).
    \end{aligned}
    \right.
\end{equation*}
A unique solution of the above system for the given ranges of $\vartheta$ and $\tau$ is
\begin{align*}
    &\tau = \pi - \arccos\frac{\ell\cos\frac{m_\ell}2 + 2\sin\frac{m_\ell}2}{\ell + m_\ell},\\
    &\vartheta = m_\ell.
\end{align*}
Here $m_\ell$ is a root of the transcendental equation\footnote{The root can be found by Newton's method with an initial guess $m_\ell^0 = \pi$, for example.}
\begin{equation*}
    (\ell + m)^2 - \left(2\sin\frac{m}2 - \ell\right)^2 - \left(2 + 2\cos\frac{m}2\right)^2 = 0, \quad m \in (0, 2\pi), \quad \ell \in (\ell_J, +\infty).
\end{equation*}

Summarizing, we obtain the critical value of angle given by
\begin{equation*}
    \vartheta_\ell^1 \overset{\mathrm{def}}{=}
    \begin{cases}
        &w_\ell, \quad \ell \in (0, \ell_J);\\
        &\vartheta_J, \quad \ell = \ell_J;\\
        &m_\ell, \quad \ell \in (\ell_J, +\infty).
    \end{cases}
\end{equation*}

In the same manner, we can determine which range of $\theta$ corresponds to the intersection $\tilde{\mathcal{B}}_\mathcal{TD}^\upsilon \cap \tilde{\mathcal{B}}_\mathcal{TD}^{ -\upsilon}$. Now we must find the critical value of the angle at which the line $\tilde{\mathcal{B}}_\mathcal{PL}^\upsilon$ crosses the surface $\tilde{\mathcal {B}}_\mathcal{TD}^{-\upsilon}$ if $\ell$ is "small". If $\ell$ is "large", then the critical value of the angle corresponds to the intersection of the universal line $\tilde{\mathcal{B}}_ \mathcal{UL}^{-\upsilon}$ with the surface $\tilde{\mathcal{B}}_\mathcal{TD}^\upsilon$. For the "medium" value of $\ell$, the critical value of the angle is $(1 - \upsilon)\pi + \upsilon\vartheta_J$. Similarly, we denote by $\vartheta^2_\ell$ a dependence of the critical value of angle from the capture radius $\ell$ for $\upsilon = +1$. If $\upsilon = -1$, this angle equals $2\pi - \vartheta^2_\ell$.

The intersection of $\tilde{\mathcal{B}}_\mathcal{PL}^\upsilon$ with $\tilde{\mathcal{B}}_\mathcal{TD}^{-\upsilon}$ is given by the system
\begin{equation*}
    \left\{
    \begin{aligned}
        &\bar{\boldsymbol{z}}_{\mathcal{B}_\mathcal{TD}^{-\upsilon}}(\tau, \vartheta) - \bar{\boldsymbol{z}}_{\mathcal{B}_\mathcal{PL}^\upsilon}(\vartheta) = \boldsymbol{0},\\
        &\vartheta \in (0, 2\pi), \quad \tau \in \left(\frac\vartheta2, \vartheta\right), \quad \ell \in (0, \ell_J).
    \end{aligned}
    \right.
\end{equation*}
Transformation $\boldsymbol{R}_\upsilon(\vartheta/2)$ of the system yields
\begin{equation*}
    \left\{
    \begin{aligned}
        &\xi_\ell(2\tau - \vartheta) - 2 - 4\cos\frac\vartheta2 = 0,\\
        &\eta_\ell(2\tau - \vartheta) + \ell = 0,\\
        &\vartheta \in (0, 2\pi), \quad \tau \in \left(\frac\vartheta2, \vartheta\right), \quad \ell \in (0, \ell_J).
    \end{aligned}
    \right.
\end{equation*}
A unique solution of the system for the given ranges of $\vartheta$ and $\tau$ is
\begin{align*}
    &\tau = \frac{w_\ell}2 + \arccos\frac{\sqrt{(\ell + w_\ell)^2 - \ell^2 + 4} - 2}4,\\
    &\vartheta = 2\arccos\frac{\sqrt{(\ell + w_\ell)^2 - \ell^2 + 4} - 2}4.
\end{align*}

The intersection of $\tilde{\mathcal{B}}_{\mathcal{UL}}^{-\upsilon}$ with $\tilde{\mathcal{B}}_\mathcal{TD}^\upsilon$ is determined by the system
\begin{equation*}
    \left\{
    \begin{aligned}
        &\bar{\boldsymbol{z}}_{\mathcal{B}_\mathcal{UL}^{-\upsilon}}(\vartheta) - \bar{\boldsymbol{z}}_{\mathcal{B}_\mathcal{TD}^\upsilon}(\tau, 2\pi - \vartheta) = \boldsymbol{0},\\
        &\vartheta \in (0, 2\pi), \quad \tau \in \left(\pi - \frac\vartheta2, 2\pi - \vartheta\right), \quad \ell \in (\ell_J, +\infty).
    \end{aligned}
    \right.
\end{equation*}
Transformation $\boldsymbol{R}_\upsilon(\vartheta/2)$ of the system yields
\begin{equation*}
    \left\{
    \begin{aligned}
        &\xi_\ell(\vartheta) - \xi_\ell(2\tau - 2\pi + \vartheta) - 4\cos\frac\vartheta2 = 0,\\ 
        &\eta_\ell(\vartheta) + \eta_\ell(2\tau - 2\pi + \vartheta) = 0,\\
        &\vartheta \in (0, 2\pi), \quad \tau \in \left(\pi - \frac\vartheta2, 2\pi - \vartheta\right), \quad \ell \in (\ell_J, +\infty).
    \end{aligned}
    \right.
\end{equation*}
A unique solution of the system for the given ranges of $\vartheta$ and $\tau$ is
\begin{align*}
    &\tau = \pi - \frac{\ell + n_\ell}2 + \frac{\sqrt{(\ell + n_\ell - 2\sin n_\ell)^2 - 4\sin^2 n_\ell}}2,\\
    &\vartheta = n_\ell.
\end{align*}
Here $n_\ell$ is a root of the transcendental equation\footnote{The root can be found by Newton's method with an initial guess $n_\ell^0 = \pi$, for example.}
\begin{equation*}
    \left\{
    \begin{aligned}
        &\eta_\ell\left(\sqrt{(\ell + n - 2\sin n)^2 - 4\sin^2n} - \ell\right) + \eta_\ell(n) = 0,\\
        &n \in (0, 2\pi), \quad \ell \in (\ell_J, +\infty).
    \end{aligned}
    \right.
\end{equation*}

Summarizing, we obtain the critical value of angle given by
\begin{equation*}
    \vartheta_\ell^2 \overset{\mathrm{def}}{=}
    \begin{cases}
        &2\arccos\frac{\sqrt{(\ell + w_\ell)^2 - \ell^2 + 4} - 2}4, \quad \ell \in (0, \ell_J);\\
        &\vartheta_J, \quad \ell = \ell_J;\\
        &n_\ell, \quad \ell \in (\ell_J, +\infty).
    \end{cases}
\end{equation*}

The critical values of angles $\vartheta^1_\ell$ and $\vartheta^2_\ell$ provide the valid ranges of $\theta$ for the computing of the intersections $\tilde{\mathcal{B}}_\mathcal{P}^\upsilon \cap \tilde{\mathcal{B}}_\mathcal{TS}^{ -\upsilon}$ and $\tilde{\mathcal{B}}_\mathcal{TD}^\upsilon \cap \tilde{\mathcal{B}}_\mathcal{TD}^{ -\upsilon}$ respectively. Moreover, for $\upsilon = +1$ an interval $(\vartheta^1_\ell, \vartheta^2_\ell)$ gives the range of $\theta$ where the surface $\tilde{\mathcal{B}}_\mathcal{P}^\upsilon$ crosses the $\tilde{\mathcal{B}}_\mathcal{TD}^{-\upsilon}$ if $\ell$ is "small" and the surface $\tilde{\mathcal{B}}_\mathcal{TD}^\upsilon$ crosses the $\tilde{\mathcal{B}}_\mathcal{TS}^{-\upsilon}$ if $\ell$ is "large" (see Fig.~\ref{fig:small_and_large_radius_slices}). For $\upsilon = -1$, the same interval is $(2\pi - \vartheta^2_\ell, 2\pi - \vartheta^1_\ell)$. Next, we will calculate these intersections.

The surface $\tilde{\mathcal{B}}_\mathcal{TS}^{-\upsilon}$ crosses the $\tilde{\mathcal{B}}_\mathcal{P}^\upsilon$ when the parameter $\vartheta \in (0, \vartheta^1_\ell)$. We use notation $\tau'$ for the parametrization of $\tilde{\mathcal{B}}_\mathcal{TS}^{-\upsilon}$ to distinguish between the parameter $\tau$ of $\tilde{\mathcal{B}}_\mathcal{P}^\upsilon$. The intersection is given by
\begin{equation*}
    \left\{
    \begin{aligned}
        &\bar{\boldsymbol{z}}_{\mathcal{B}_\mathcal{TS}^{-\upsilon}}\left(\tau', \vartheta\right) - \bar{\boldsymbol{z}}_{\mathcal{B}_\mathcal{P}^\upsilon}\left(\tau, \vartheta\right) = \boldsymbol{0},\\
        &\vartheta \in (0, \vartheta^1_\ell), \quad \tau \in \left(0, \pi - \frac\vartheta2\right), \quad \tau' \in (\vartheta, 2\pi).
    \end{aligned}
    \right.
\end{equation*}
Transformation $\boldsymbol{R}_\upsilon(\vartheta/2)$ of the system leads to 
\begin{equation*}
    \left\{
    \begin{aligned}
        &(\ell + \vartheta)\sin\left(\tau' - \frac\vartheta2\right) - 2\cos\frac\vartheta2 + 2\cos\left(\tau + \frac\vartheta2\right) = 0,\\
        &(\ell + \vartheta)\cos\left(\tau' - \frac\vartheta2\right) - 2\sin\frac\vartheta2 + \ell = 0,\\
        &\vartheta \in (0, \vartheta^1_\ell), \quad \tau \in \left(0, \pi - \frac\vartheta2\right), \quad \tau' \in (\vartheta, 2\pi).
    \end{aligned}
    \right.
\end{equation*}
For $\vartheta \in (0, \vartheta^1_\ell)$, a unique solution of the system for the given ranges of $\tau$ and $\tau'$ is
\begin{equation}\label{eq:BP_BTS_intersection}
    \begin{split}
        &\tau = -\frac\vartheta2 + \arccos\left(\cos\frac\vartheta2 - \frac12\sqrt{(\ell+\vartheta)^2 - \left(\ell-2\sin\frac\vartheta2\right)^2}\right),\\
        &\tau' = \frac\vartheta2 + \arccos\frac{2\sin\frac\vartheta2 - \ell}{\ell + \vartheta}.
    \end{split}
\end{equation}

Similarly, the intersection of $\tilde{\mathcal{B}}_\mathcal{TD}^\upsilon$ and $\tilde{\mathcal{B}}_\mathcal{TD}^{-\upsilon}$ is given by the system
\begin{equation*}
    \left\{
    \begin{aligned}
        &\bar{\boldsymbol{z}}_{\mathcal{B}_\mathcal{TD}^{-\upsilon}}\left(\tau, \vartheta\right) - \bar{\boldsymbol{z}}_{\mathcal{B}_\mathcal{TD}^\upsilon}\left(\tau', 2\pi - \vartheta\right) = \boldsymbol{0},\\
        &\vartheta \in (\vartheta^2_\ell, 2\pi - \vartheta^2_\ell), \quad \tau \in \left(\frac\vartheta2, \vartheta\right), \quad \tau' \in \left(\pi - \frac\vartheta2, 2\pi - \vartheta\right).\\
    \end{aligned}
    \right.
\end{equation*}
The parameter $\tau$ relates to $\tilde{\mathcal{B}}_\mathcal{TD}^{-\upsilon}$ and $\tau'$ to $\tilde{\mathcal{B}}_\mathcal{TD}^\upsilon$. Transformation $\boldsymbol{R}_\upsilon(\vartheta/2)$ of the system yields
\begin{equation*}
    \left\{
    \begin{aligned}
        &\xi_\ell(2\tau - \vartheta) - \xi_\ell(2\tau' + \vartheta - 2\pi) - 4\cos\frac\vartheta2 = 0,\\
        &\eta_\ell(2\tau - \vartheta) + \eta_\ell(2\tau' + \vartheta - 2\pi) = 0,\\
        &\vartheta \in (\vartheta^2_\ell, 2\pi - \vartheta^2_\ell), \quad \tau \in \left(\frac\vartheta2, \vartheta\right), \quad \tau' \in \left(\pi - \frac\vartheta2, 2\pi - \vartheta\right).\\
    \end{aligned}
    \right.
\end{equation*}
For $\vartheta \in (\vartheta^2_\ell, 2\pi - \vartheta^2_\ell)$, a unique solution of the system for the given ranges of $\tau$ and $\tau'$ is
\begin{equation}\label{eq:BTD_BTD_intersection}
    \begin{split}
        &\tau = \frac{\vartheta + p_\ell(\vartheta)}2,\\
        &\tau' = \pi - \frac{\ell + \vartheta}2 + \frac12\sqrt{\left(\xi_\ell(p_\ell(\vartheta)) - 4\cos\frac\vartheta2\right)^2 + \eta_\ell^2(p_\ell(\vartheta)) - 4},
    \end{split}
\end{equation}
where $p_\ell(\vartheta)$ is a root of the transcendental equation\footnote{The root can be found by Newton's method with an initial guess $p^0_\ell(\vartheta) = \vartheta$, for example.}
\begin{equation*}
    \left\{
    \begin{aligned}
        &\eta_\ell(p) + \eta_\ell\left(\sqrt{\left(\xi_\ell(p) - 4\cos\frac\vartheta2\right)^2 + \eta_\ell^2(p) - 4} - \ell\right) = 0,\\
        &p \in (0, \vartheta), \quad \vartheta \in (\vartheta^2_\ell, 2\pi - \vartheta^2_\ell).
    \end{aligned}
    \right.
\end{equation*}

As noted previously, the surface $\tilde{\mathcal{B}}_\mathcal{TD}^{-\upsilon}$ crosses the $\tilde{\mathcal{B}}_\mathcal{P}^\upsilon$ only if $\ell \in (0, \ell_J)$. We associate a parameter $\tau$ with $\tilde{\mathcal{B}}_\mathcal{P}^\upsilon$ and $\tau'$ with $\tilde{\mathcal{B}}_\mathcal{TD}^{-\upsilon}$. Consider the system
\begin{equation*}
    \left\{
    \begin{aligned}
        &\bar{\boldsymbol{z}}_{\mathcal{B}_\mathcal{TD}^{-\upsilon}}\left(\tau', \vartheta\right) - \bar{\boldsymbol{z}}_{\mathcal{B}_\mathcal{P}^\upsilon}\left(\tau, \vartheta\right) = \boldsymbol{0},\\
        &\vartheta \in (\vartheta^1_\ell, \vartheta^2_\ell), \quad \tau \in \left(0, \pi - \frac\vartheta2\right), \quad \tau' \in \left(\frac\vartheta2, \vartheta\right), \quad \ell \in (0, \ell_J).\\
    \end{aligned}
    \right.
\end{equation*}
Transformation $\boldsymbol{R}_\upsilon(\vartheta/2)$ of the system leads to 
\begin{equation*}
    \left\{
    \begin{aligned}
        &\xi_\ell(2\tau' - \vartheta) + 2\cos\frac{2\tau + \vartheta}2 - 4\cos\frac\vartheta2 = 0,\\
        &\eta_\ell(2\tau' - \vartheta) + \ell = 0,\\
        &\vartheta \in (\vartheta^1_\ell, \vartheta^2_\ell), \quad \tau \in \left(0, \pi - \frac\vartheta2\right), \quad \tau' \in \left(\frac\vartheta2, \vartheta\right), \quad \ell \in (0, \ell_J).\\
    \end{aligned}
    \right.
\end{equation*}
For $\vartheta \in (\vartheta^1_\ell, \vartheta^2_\ell)$, a unique solution of the system for the given ranges of $\tau$ and $\tau'$ is
\begin{equation}\label{eq:BP_BTD_intersection}
    \begin{aligned}
        &\tau = -\frac\vartheta2 + \arccos\left(2\cos\frac\vartheta2 - \frac{\sqrt{(\ell + w_\ell)^2 - \ell^2 + 4}}2\right),\\
        &\tau' = \frac{w_\ell + \vartheta}2.
    \end{aligned}
\end{equation}

The surface $\tilde{\mathcal{B}}_\mathcal{TD}^\upsilon$ crosses the $\tilde{\mathcal{B}}_\mathcal{TS}^{-\upsilon}$ only if $\ell \in (\ell_J, +\infty)$. As earlier, we associate parameters $\tau$, $\tau'$ with $\tilde{\mathcal{B}}_\mathcal{TS}^{-\upsilon}$, $\tilde{\mathcal{B}}_\mathcal{TD}^\upsilon$ respectively. Consider the system 
\begin{equation*}
    \left\{
    \begin{aligned}
        &\bar{\boldsymbol{z}}_{\mathcal{B}_\mathcal{TS}^{-\upsilon}}\left(\tau, \vartheta\right) - \bar{\boldsymbol{z}}_{\mathcal{B}_\mathcal{TD}^\upsilon}\left(\tau', 2\pi - \vartheta\right) = \boldsymbol{0},\\
        &\vartheta \in (\vartheta^1_\ell, \vartheta^2_\ell), \quad \tau \in (\vartheta, 2\pi), \quad \tau' \in \left(\pi - \frac\vartheta2, 2\pi - \vartheta\right), \quad \ell \in (\ell_J, +\infty).
    \end{aligned}
    \right.
\end{equation*}
Transformation $\boldsymbol{R}_\upsilon(\vartheta/2)$ of the system gives
\begin{equation*}
    \left\{
    \begin{aligned}
        &(\ell + \vartheta)\sin\left(\tau - \frac\vartheta2\right) - 2\cos\frac\vartheta2 - \xi_\ell(2\tau' + \vartheta - 2\pi) = 0,\\
        &(\ell + \vartheta)\cos\left(\tau - \frac\vartheta2\right) - 2\sin\frac\vartheta2 + \eta_\ell(2\tau' + \vartheta - 2\pi) = 0,\\
        &\vartheta \in (\vartheta^1_\ell, \vartheta^2_\ell), \quad \tau \in (\vartheta, 2\pi), \quad \tau' \in \left(\pi - \frac\vartheta2, 2\pi - \vartheta\right), \quad \ell \in (\ell_J, +\infty).
    \end{aligned}
    \right.
\end{equation*}
For $\vartheta \in (\vartheta^1_\ell, \vartheta^2_\ell)$, a unique solution of the system for the given ranges of $\tau$ and $\tau'$ is given by
\begin{equation}\label{eq:BTS_BTD_intersection}
    \begin{aligned}
        &\tau = \pi - \arcsin\frac{(\ell + \vartheta)^2 - (\ell + q_\ell(\vartheta))^2}{4(\ell + \vartheta)},\\
        &\tau' = \pi + \frac{q_\ell(\vartheta) - \vartheta}2, 
    \end{aligned}
\end{equation}
where $q_\ell(\vartheta)$ is a root of the transcendental equation\footnote{The root can be found by Newton's method with an initial guess $q^0_\ell(\vartheta) = \vartheta$, for example.}
\begin{equation}\label{eq:q_ell}
    \left\{
    \begin{aligned}
        &(\ell + \vartheta)^2 - \left(2 + 2\cos\frac{q - \vartheta}2\right)^2 - \left(\ell + q + 2\sin\frac{q - \vartheta}2\right)^2 = 0,\\
        &q \in (0, 2\pi - \vartheta), \quad \vartheta \in (\vartheta^1_\ell, \vartheta^2_\ell), \quad \ell \in (\ell_J, +\infty).
    \end{aligned}
    \right.
\end{equation}

\subsection{Actual parts of the barrier}

Now when all conditions of intersection are obtained we can determine the valid parts of the emanated semipermeable surfaces that constitute barrier $\mathcal{B}$. At first, we describe surface $\mathcal{B}_\mathcal{P}^\upsilon = \tilde{\mathcal{B}}_\mathcal{P}^\upsilon \cap \mathcal{B}$. The $\theta$-slices of this surface exist for all $\vartheta \in (0, 2\pi)$ (see Figs. \ref{fig:small_and_large_radius_slices}-\ref{fig:middle_radius_slices}). For all values of $\ell$, if $\vartheta \in (0, \vartheta^1_\ell]$, the intersection is determined by \eqref{eq:BP_BTS_intersection}. If $\ell \in (0, \ell_J)$ the surface $\tilde{\mathcal{B}}_\mathcal{P}^\upsilon$ also crosses $\tilde{\mathcal{B}}_\mathcal{TD}^{-\upsilon}$ and the maximal value of $\tau$ is given by \eqref{eq:BP_BTD_intersection}. Summarizing, we have
\begin{equation}\label{eq:BPv_tau_vartheta}
    \mathcal{B}_\mathcal{P}^\upsilon \overset{\mathrm{def}}{=} \left\{\bar{\boldsymbol{z}}_{\mathcal{B}_\mathcal{P}^\upsilon}(\tau, \vartheta):\: \vartheta \in (0, 2\pi), \quad \tau \in \left(0, \tau^{\max}_{\mathcal{B}_\mathcal{P}, \ell}(\vartheta)\right)\right\},
\end{equation}
where
\begin{equation*}
    \tau^{\max}_{\mathcal{B}_\mathcal{P}, \ell}(\vartheta) \overset{\mathrm{def}}{=}
    \begin{cases}
        &\arccos\left(\cos\frac\vartheta2 - \frac{\sqrt{(\ell + \vartheta)^2 + (\ell - 2\sin\frac\vartheta2)^2}}2\right) - \frac\vartheta2, \quad \vartheta \in \left(0, \vartheta^1_\ell\right];\\
        &\arccos\left(2\cos\frac\vartheta2 - \frac{\sqrt{(\ell + w_\ell)^2 - \ell^2 + 4}}2\right) - \frac\vartheta2, \quad \vartheta \in \left[\vartheta^1_\ell, \vartheta^{21}_\ell\right];\\
        &\pi - \frac\vartheta2, \quad \vartheta \in \left[\vartheta^{21}_\ell, 2\pi\right)
    \end{cases}
\end{equation*}
and
\begin{equation*}
    \vartheta^{21}_\ell \overset{\mathrm{def}}{=}
    \begin{cases}
        &\vartheta^2_\ell \quad \ell \in (0, \ell_J];\\
        &\vartheta^1_\ell, \quad \ell \in [\ell_J, +\infty).\\
    \end{cases}
\end{equation*}

Next, we will describe the surface $\mathcal{B}_\mathcal{TS}^\upsilon = \tilde{\mathcal{B}}_\mathcal{TS}^\upsilon \cap \mathcal{B}$. The $\theta$-slices of this surface exist only for $\vartheta \in (0, \vartheta^1_\ell)$ if $\ell \in (0, \ell_J]$ and for $\vartheta \in (0, \vartheta^2_\ell)$ if $\ell \in (\ell_J, +\infty)$ (see Figs. \ref{fig:small_and_large_radius_slices}-\ref{fig:middle_radius_slices}). The maximal value of time-to-go $\tau$ is determined by $\tau'$ of \eqref{eq:BP_BTS_intersection} and $\tau$ of \eqref{eq:BTS_BTD_intersection}. Thus,
\begin{equation}\label{eq:BTSv_tau_vartheta}
    \mathcal{B}_\mathcal{TS}^\upsilon \overset{\mathrm{def}}{=} \left\{\bar{\boldsymbol{z}}_{\mathcal{B}_\mathcal{TS}^\upsilon}(\tau, \vartheta):\: \vartheta \in (0, \vartheta^{12}_\ell), \quad \tau \in \left(\vartheta, \tau_{\mathcal{B}_\mathcal{TS}, \ell}^{\max}(\vartheta)\right)\right\},
\end{equation}
where
\begin{equation*}
    \vartheta^{12}_\ell \overset{\mathrm{def}}{=}
    \begin{cases}
        &\vartheta^1_\ell, \quad \ell \in (0, \ell_J];\\
        &\vartheta^2_\ell, \quad \ell \in [\ell_J, +\infty)
    \end{cases}
\end{equation*}
and
\begin{equation*}
    \tau^{\max}_{\mathcal{B}_\mathcal{TS}, \ell}(\vartheta) \overset{\mathrm{def}}{=}
    \begin{cases}
        &\frac{\vartheta}2 + \arccos\frac{2\sin\frac\vartheta2 - \ell}{\ell + \vartheta}, \quad \vartheta \in (0, \vartheta^1_\ell];\\
        &\pi - \arcsin\frac{(\ell + \vartheta)^2 - (\ell + q_\ell(\vartheta))^2}{4(\ell + \vartheta)}, \quad \vartheta \in [\vartheta^1_\ell, \vartheta^{12}_\ell].
    \end{cases}
\end{equation*}
\begin{lemma}\label{lem:tau_max_BTD}
    The next interval estimation is valid: $0 < \tau^{\max}_{\mathcal{B}_\mathcal{TS}, \ell}(\vartheta) < \vartheta + \pi$ for all $\vartheta \in (0, \vartheta^{12}_\ell]$.
\end{lemma}
\begin{proof}
   For $\vartheta \in (0, \vartheta^1_\ell]$ (or $\ell \leq \ell_J)$, the statement is trivial since
    \begin{equation*}
        \tau^{\max}_{\mathcal{B}_\mathcal{TS}, \ell}(\vartheta) - \vartheta = -\frac\vartheta2 + \arccos\frac{2\sin\frac\vartheta2 - \ell}{\ell + \vartheta} < \arccos\frac{2\sin\frac\vartheta2 - \ell}{\ell + \vartheta} \leq \pi.
    \end{equation*}
    For $\vartheta \in (\vartheta^1_\ell, \vartheta^2_\ell]$ (and $\ell > \ell_J)$, it is suffice to prove that $q_\ell(\vartheta) < \vartheta$ since
    \begin{equation*}
        \tau^{\max}_{\mathcal{B}_\mathcal{TS}, \ell}(\vartheta) - \vartheta < \tau^{\max}_{\mathcal{B}_\mathcal{TS}, \ell}(\vartheta) = \pi - \arcsin\frac{(\vartheta - q_\ell(\vartheta))(2\ell + \vartheta + q_\ell(\vartheta))}{4(\ell + \vartheta)} < \pi.
    \end{equation*}
    Assume the contrary that $q_\ell(\vartheta) \geq \vartheta$. Transform \eqref{eq:q_ell} to
    \begin{equation*}
        (\ell + \vartheta)^2 - (\ell + q_\ell(\vartheta))^2 = 8\left(1 + \cos\frac{q_\ell(\vartheta) - \vartheta}2\right) + 4(\ell + q_\ell(\vartheta))\sin\frac{q_\ell(\vartheta) - \vartheta}2.
    \end{equation*}
    According to \eqref{eq:q_ell}, the root $q_\ell(\vartheta)$ is less than $2\pi - \vartheta$ which implies $0 \leq (q_\ell(\vartheta) - \vartheta)/2 < \pi$. Hence, the left part of the equation is non-positive while the right is positive. The contradiction completes the proof. \qed
\end{proof}

Further, we will obtain a description of $\mathcal{B}_\mathcal{TD}^\upsilon = \tilde{\mathcal{B}}_\mathcal{TD}^\upsilon \cap \mathcal{B}$. The $\theta$-slices of this surface exist only for $\vartheta \in (0, 2\pi - \vartheta^2_\ell)$ if $\ell \in (0, \ell_J]$ and for $\vartheta \in (0, 2\pi - \vartheta^1_\ell)$ if $\ell \in [\ell_J, +\infty)$ (see Figs. \ref{fig:small_and_large_radius_slices}-\ref{fig:middle_radius_slices}). The maximal value of time-to-go $\tau$ is determined by $\tau'$ of \eqref{eq:BP_BTD_intersection}, $\tau'$ of \eqref{eq:BTS_BTD_intersection}, and $\tau$, $\tau'$ of \eqref{eq:BTD_BTD_intersection}. Therefore,
\begin{equation}\label{eq:BTDv_tau_vartheta}
    \mathcal{B}_\mathcal{TD}^\upsilon \overset{\mathrm{def}}{=} \left\{\bar{\boldsymbol{z}}_{\mathcal{B}_\mathcal{TD}^\upsilon}(\tau, \vartheta):\: \vartheta \in (0, 2\pi - \vartheta^{21}_\ell), \quad \tau \in \left(\frac\vartheta2, \tau_{\mathcal{B}_\mathcal{TD}, \ell}^{\max}(\vartheta)\right)\right\},
\end{equation}
where
\begin{equation*}
    \tau^{\max}_{\mathcal{B}_\mathcal{TD}, \ell}(\vartheta) \overset{\mathrm{def}}{=}
    \begin{cases}
        &\vartheta, \quad \vartheta \in (0, \vartheta^{12}_\ell];\\
        &\frac{w_\ell + \vartheta}2, \quad \vartheta \in [\vartheta^1_\ell, \vartheta^2_\ell], \quad \ell \in (0, \ell_J];\\
        &\frac{p_\ell(\vartheta) + \vartheta}2, \quad \vartheta \in [\vartheta^2_\ell, 2\pi - \vartheta^2_\ell];\\
        &\frac{q_\ell(2\pi - \vartheta) + \vartheta}2, \quad \vartheta \in [2\pi - \vartheta^2_\ell, 2\pi - \vartheta^1_\ell], \quad \ell \in [\ell_J, +\infty).\\
    \end{cases}
\end{equation*}

Finally, we give expressions for $\mathcal{B}_\mathcal{PL}^\upsilon = \tilde{\mathcal{B}}_\mathcal{PL}^\upsilon \cap \mathcal{B}$ and $\mathcal{B}_\mathcal{UL}^\upsilon = \tilde{\mathcal{B}}_\mathcal{UL}^\upsilon \cap \mathcal{B}$ using the fact that these lines are parts of the boundaries of the corresponding surfaces:
\begin{equation}\label{eq:BPL_BUL_vartheta}
    \mathcal{B}_\mathcal{PL}^\upsilon \overset{\mathrm{def}}{=} \left\{\bar{\boldsymbol{z}}_{\mathcal{B}_\mathcal{PL}^\upsilon}(\vartheta):\: \vartheta \in (\vartheta^{21}_\ell, 2\pi)\right\}, \: \mathcal{B}_\mathcal{UL}^\upsilon \overset{\mathrm{def}}{=} \left\{\bar{\boldsymbol{z}}_{\mathcal{B}_\mathcal{UL}^\upsilon}(\vartheta):\: \vartheta \in (0, \vartheta^{12}_\ell]\right\}.
\end{equation}
For the sake of completeness, we also provide a parametric description of the barrier dispersal line $\mathcal{B}_\mathcal{DL}$. In the description, we use the function
\begin{equation*}
    \gamma_\theta \overset{\mathrm{def}}{=}
    \begin{cases}
        &-1, \quad \theta \in (0, \pi];\\
        &+1, \quad \theta \in [\pi, 2\pi).
    \end{cases}
\end{equation*}
For all values of capture radius $\ell$ and all $\theta$-slices, the dispersal line is a part of a boundary of $\mathcal{B}_\mathcal{TS}^\upsilon$ or $\mathcal{B}_\mathcal{TD}^\upsilon$. Using the fact, we obtain
\begin{equation}\label{eq:BDL_theta}
    \mathcal{B}_\mathcal{DL} \overset{\mathrm{def}}{=} \left\{\boldsymbol{z}_{\mathcal{B}_\mathcal{DL}}(\theta):\: \theta \in (0, 2\pi)\right\},
\end{equation}
where
\begin{equation*}
    \boldsymbol{z}_{\mathcal{B}_\mathcal{DL}}(\theta) \overset{\mathrm{def}}{=}
    \begin{cases}
        &\bar{\boldsymbol{z}}_{\mathcal{B}_\mathcal{TD}^{\gamma_\theta}}\left(\tau^{\max}_{\mathcal{B}_\mathcal{TD}}\left(\pi - |\pi - \theta|\right), \pi - |\pi - \theta|\right), \quad \theta \in [\vartheta_\ell^{12}, 2\pi - \vartheta_\ell^{12}];\\
        &\bar{\boldsymbol{z}}_{\mathcal{B}_\mathcal{TS}^{\gamma_\theta}}\left(\tau^{\max}_{\mathcal{B}_\mathcal{TS}}\left(\pi - |\pi - \theta|\right), \pi - |\pi - \theta|\right), \quad \text{otherwise}.
    \end{cases}
\end{equation*}

\begin{figure}
    \begin{center}
        \includegraphics[width=0.44\linewidth]{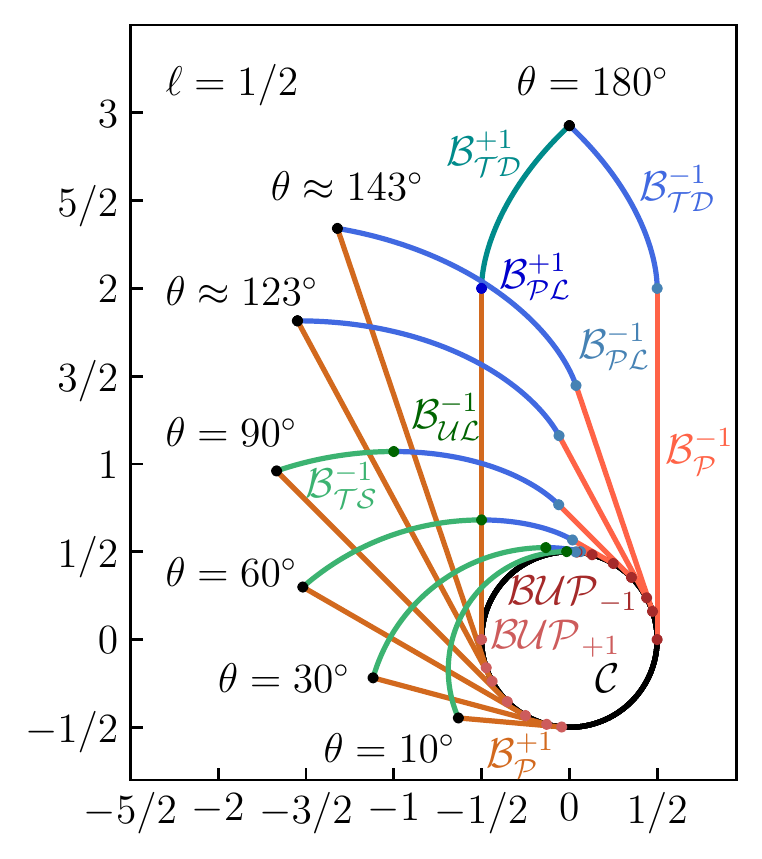}
        \includegraphics[width=0.47\linewidth]{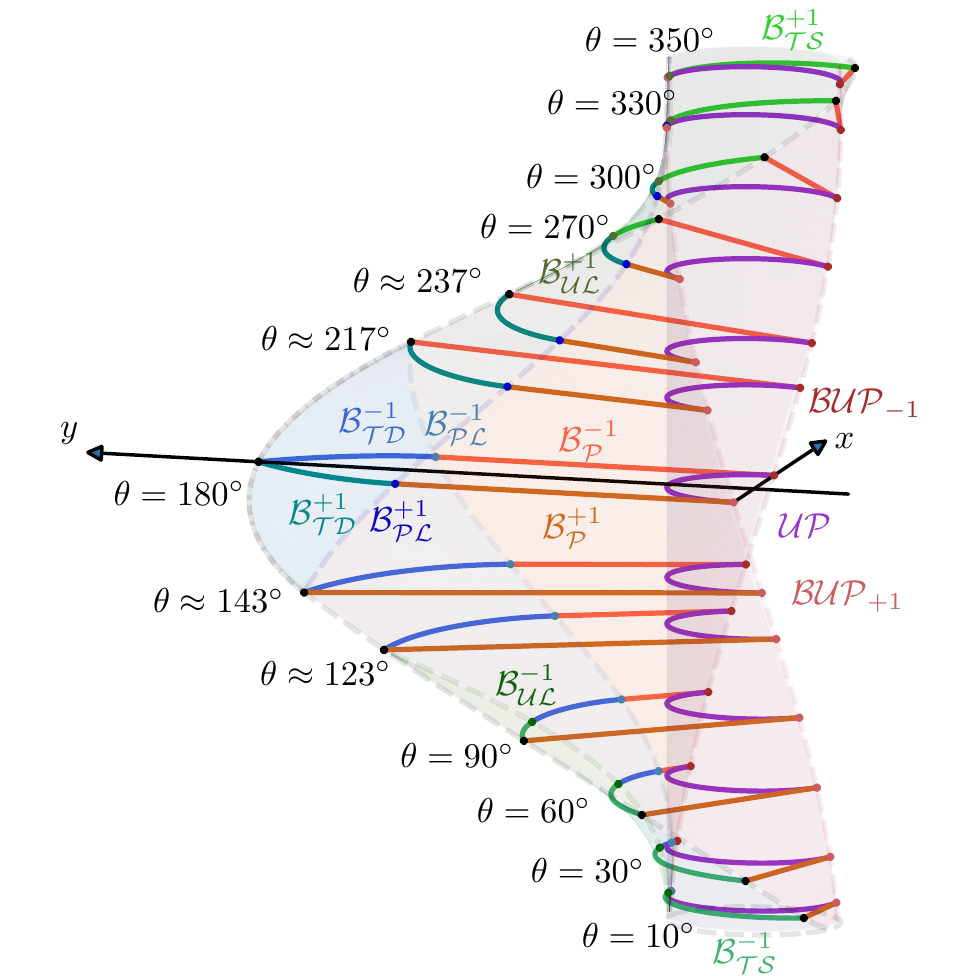}
        \includegraphics[width=0.44\linewidth]{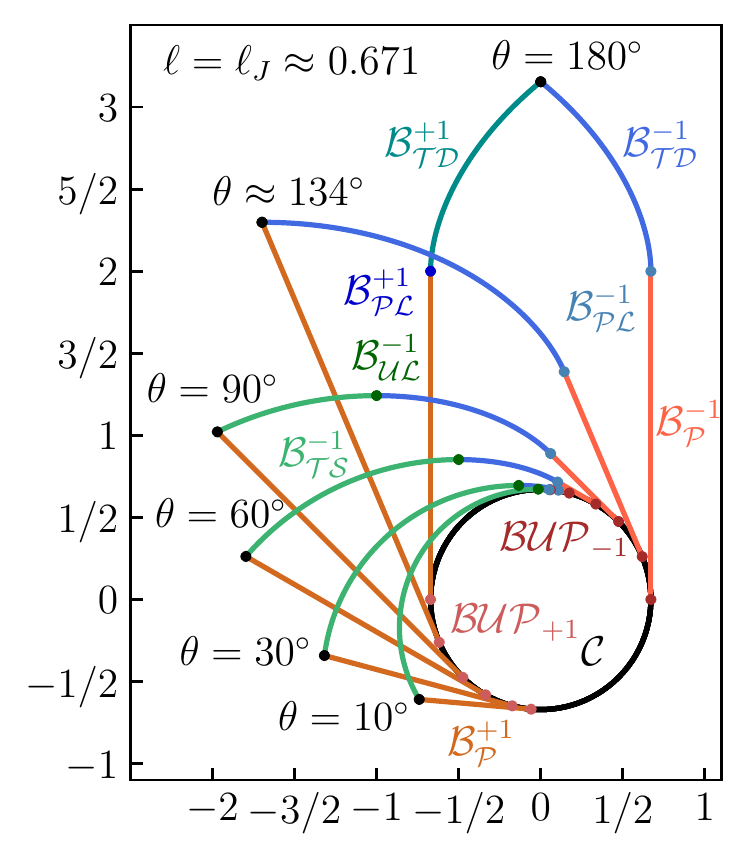}
        \includegraphics[width=0.47\linewidth]{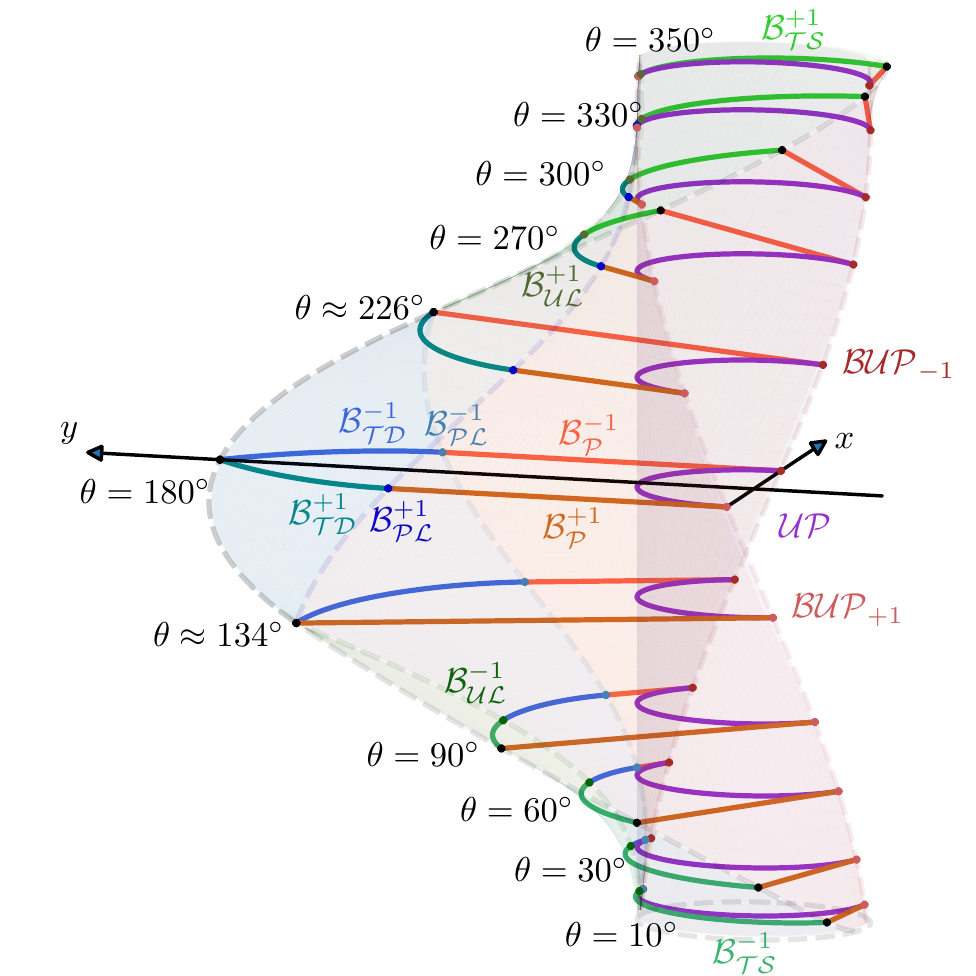}
        \includegraphics[width=0.44\linewidth]{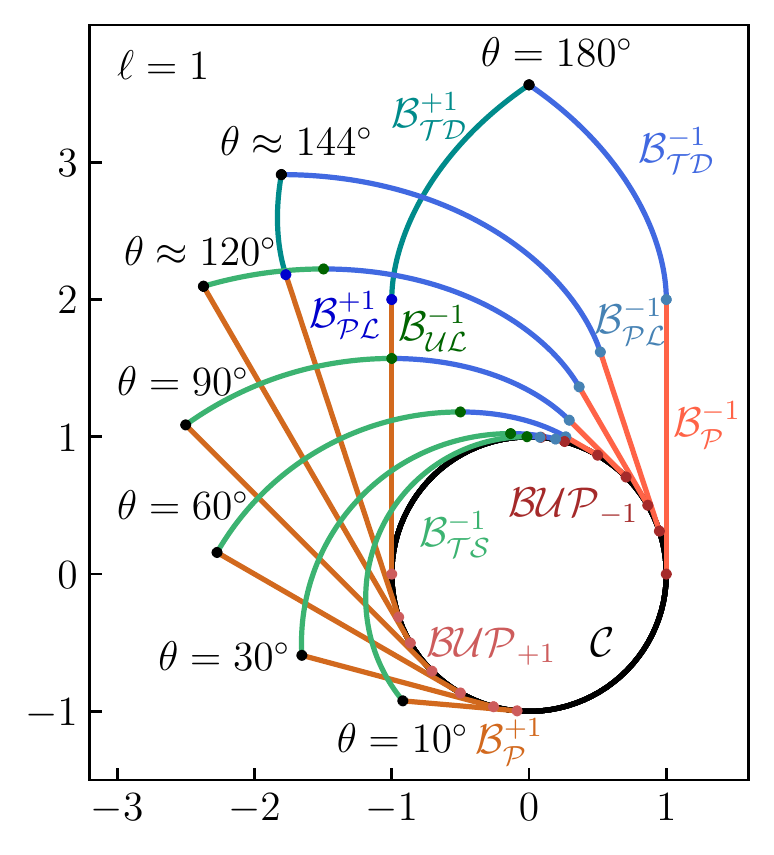}
        \includegraphics[width=0.47\linewidth]{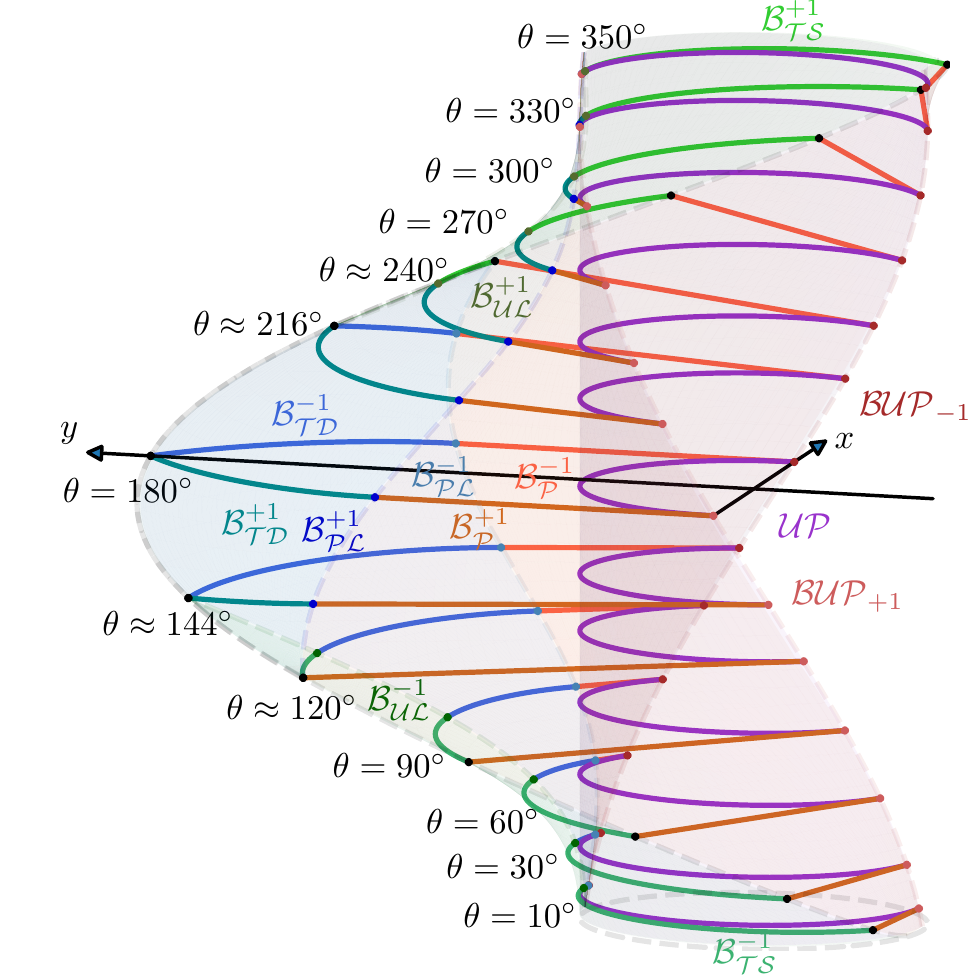}
        \caption{Barrier cross sections for the "small", "medium", and "large" capture radii}
        \label{fig:all_slices}
    \end{center}
\end{figure}
\begin{figure}
    \begin{center}
        \includegraphics[width=0.31\linewidth]{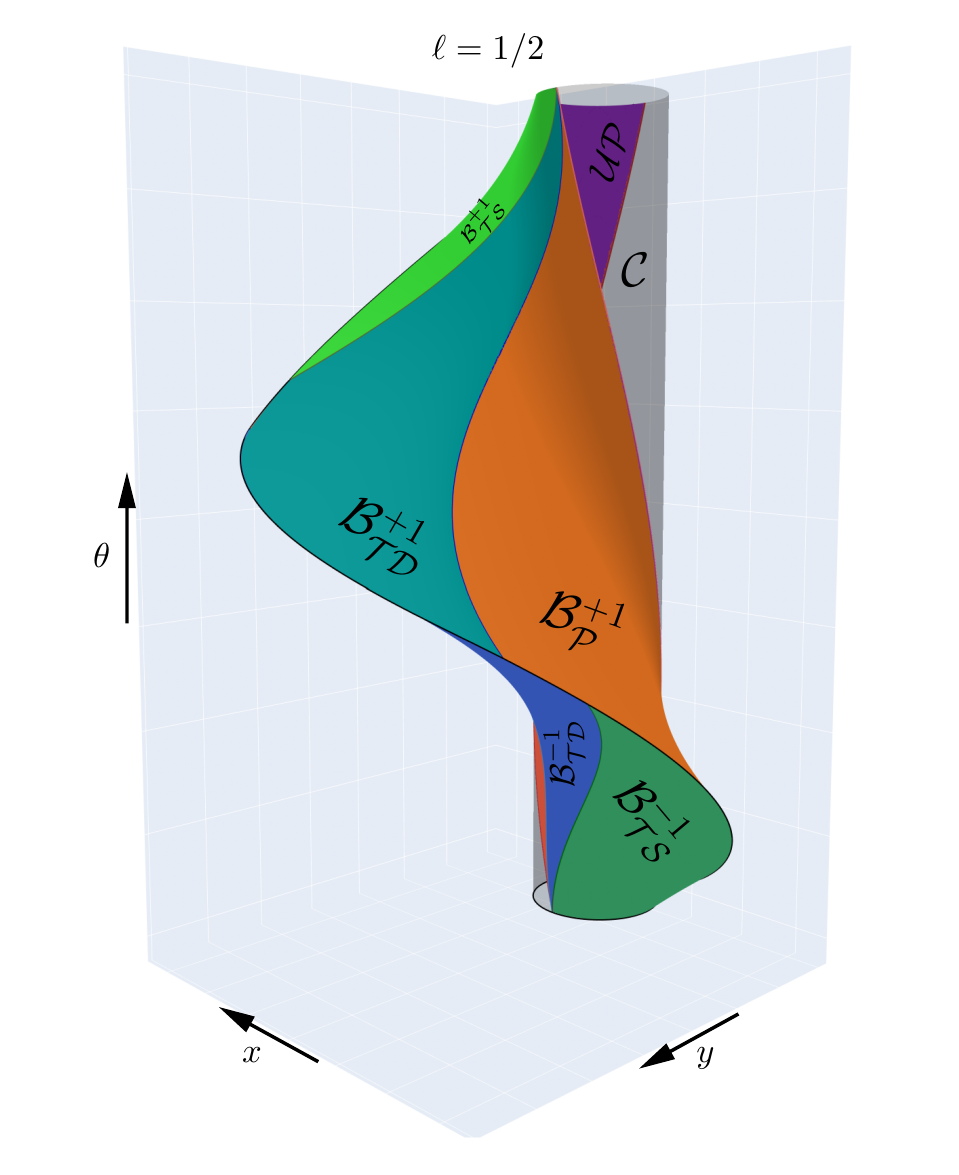}
        \includegraphics[width=0.32\linewidth]{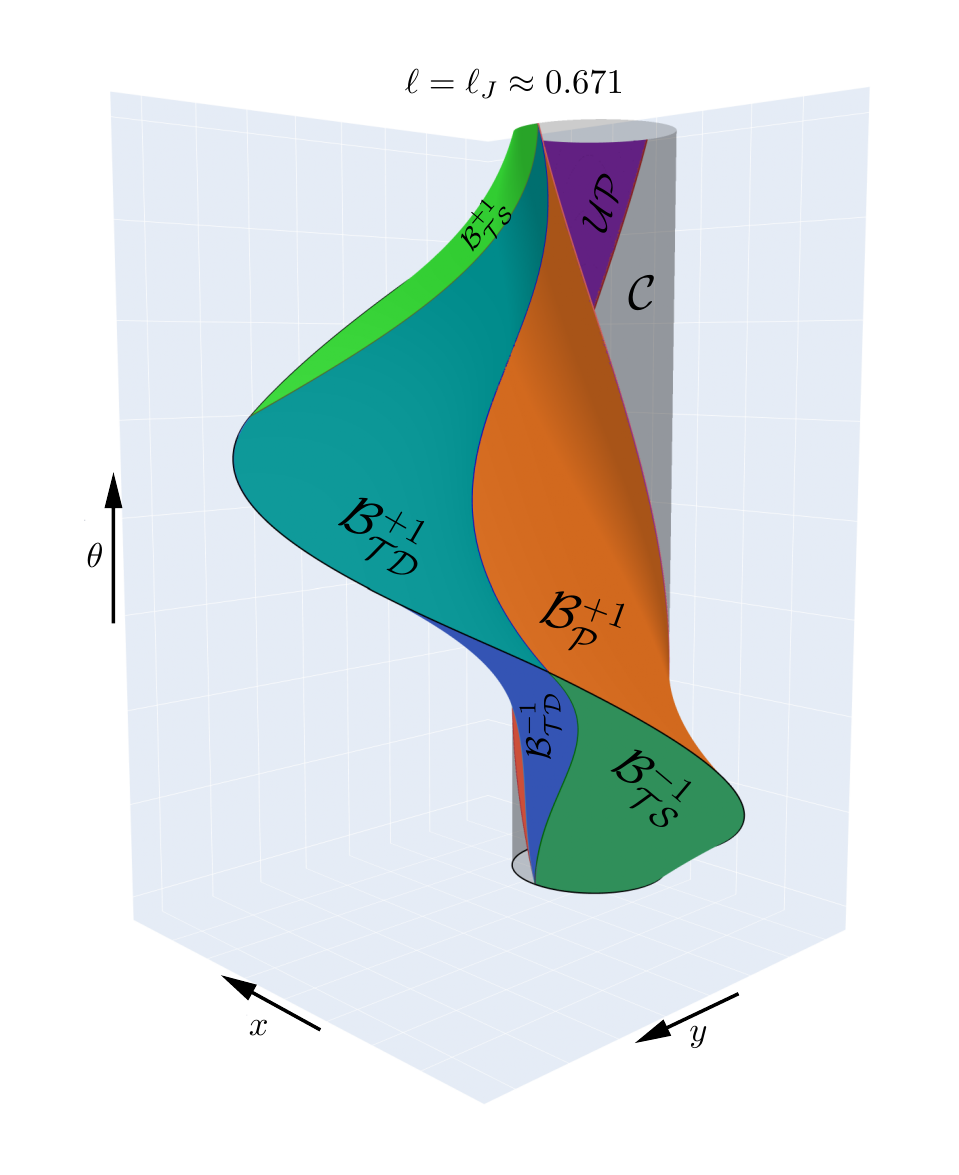}
        \includegraphics[width=0.33\linewidth]{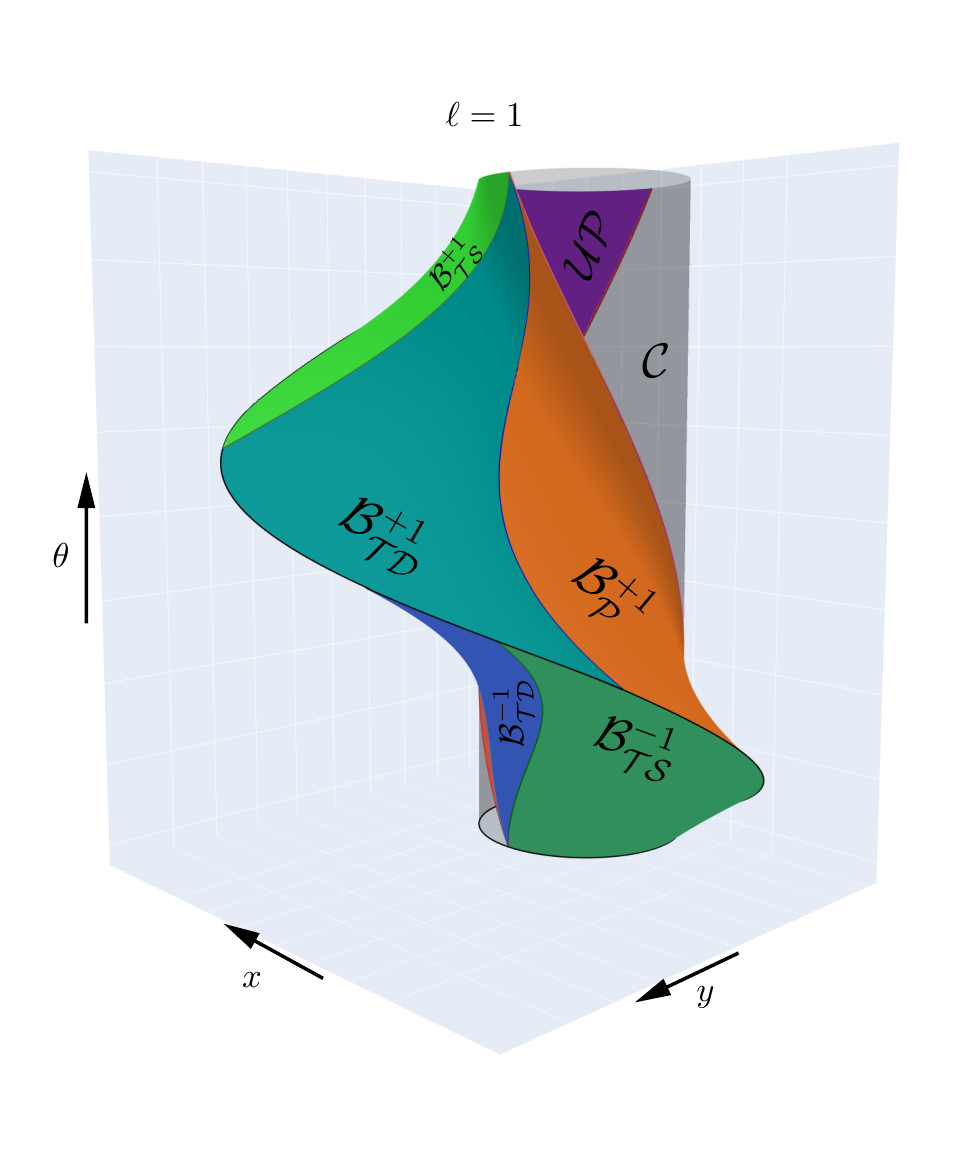}
        \caption{Difference between the "small" ($\ell < \ell_J$), "medium" ($\ell = \ell_J$), and "large" ($\ell > \ell_J$) capture radii cases}
        \label{fig:all}
    \end{center}
\end{figure}

Visualization of the surfaces $\mathcal{B}_\mathcal{P}^\upsilon$, $\mathcal{B}_\mathcal{TS}^\upsilon$, $\mathcal{B}_\mathcal{TD}^\upsilon$ and lines $\mathcal{B}_\mathcal{PL}^\upsilon$, $\mathcal{B}_\mathcal{UL}^\upsilon$ for the different values of the capture radius $\ell$ is presented in Figs~\ref{fig:all_slices}. These figures demonstrate $\theta$-slices of the 3-dimensional barrier surface $\mathcal{B}$. Colors of the lines and dots correspond to colored labels with the names of surfaces and lines. In Fig.~\ref{fig:all} we highlight a difference between the "small", "medium", and "large" radius cases. If the capture radius is "small", then $\mathcal{B}_\mathcal{P}^{+1}$ and $\mathcal{B}_\mathcal{TD}^{-1}$ have a union part of their boundary. For the "large" capture radius, the union part is between $\mathcal{B}_\mathcal{TS}^{-1}$ and $\mathcal{B}_\mathcal{TD}^{+1}$. For the "medium" capture radius, the union part degenerates into a point where $\mathcal{B}_\mathcal{P}^{+1}$, $\mathcal{B}_\mathcal{TD}^{-1}$, $\mathcal{B}_\mathcal{TS}^{-1}$, and $\mathcal{B}_\mathcal{TD}^{+1}$ meet each other.

\section{Optimal feedback controls on the barrier}\label{sec:optimal_feedback_controls}

In the previous sections, we have obtained parametric descriptions of all pieces of the barrier with corresponding restrictions on parameters. Since the values of optimal controls for both players on each of these pieces are known, by checking which of the pieces the current state belongs to, we can calculate the values of optimal controls for both players. For the obtained parametrizations \eqref{eq:BPv_tau_vartheta}-\eqref{eq:BPL_BUL_vartheta}, such verification is complicated by additional calculation of the parameters $\tau$ and $\vartheta$. Only for the dispersal line $\mathcal{B}_{\mathcal{DL}}$, the parameterization \eqref{eq:BDL_theta} immediately allows checking whether the state $\boldsymbol{z} = \begin{bmatrix} x & y & \theta \end{bmatrix}$ belongs to the dispersal line. In this section, we will obtain a parametrization for \eqref{eq:BPv_tau_vartheta}-\eqref{eq:BPL_BUL_vartheta} via the state vector $\boldsymbol{z}$.

For the GTIC it will be demonstrated that each piece of the barrier (we denote the piece by $\mathcal{P}$) can be parameterized in the following form
\begin{equation}\label{eq:special_form_of_pieces}
    \mathcal{P} = \{\boldsymbol{z} \in \mathcal{F}^{\ell}_\mathcal{P}:\: \ell = \ell_\mathcal{P}(\boldsymbol{z})\}.
\end{equation}
Here $\mathcal{F}_\mathcal{P}^\ell \subset \mathbb{R}^2 \times \mathbb{S}$ is a frame set. It corresponds to restrictions on parameters $\tau$ and $\vartheta$.

To obtain the desired parametric description of $\mathcal{B}_\mathcal{P}^{\upsilon}$ we have to exclude $\tau$ and $\vartheta$ from the system of equations
\begin{equation*}
    \left\{
    \begin{aligned}
        &\bar{\boldsymbol{z}}_{\mathcal{B}_\mathcal{P}^{\upsilon}}(\tau, \vartheta) = \boldsymbol{z};\\
        &\vartheta \in (0, 2\pi), \: \tau \in \left(0, \tau^{\max}_{\mathcal{B}_\mathcal{P}, \ell}(\vartheta)\right).
    \end{aligned}
    \right.
\end{equation*}
Applying the transformation $\boldsymbol{R}_\upsilon(\vartheta/2)$ we have
\begin{equation*}
    \left\{
    \begin{aligned}
        &2\cos\frac\vartheta2 - 2\cos\left(\tau + \frac\vartheta2\right) = -\upsilon x\cos\frac\vartheta2 + y\sin\frac\vartheta2;\\
        &-\ell = \upsilon x\sin\frac\vartheta2 + y\sin\frac\vartheta2;\\
        &(1 - \upsilon)\pi + \upsilon\vartheta = \theta; \quad \vartheta \in (0, 2\pi), \: \tau \in \left(0, \tau^{\max}_{\mathcal{B}_\mathcal{P}, \ell}(\vartheta)\right).
    \end{aligned}
    \right.
\end{equation*}
Eliminating $\vartheta$ from the second equation we obtain
\begin{equation*}
    \ell = -\upsilon\left(x\sin\frac\theta2 + y\cos\frac{\theta}2\right) \overset{\mathrm{def}}{=} \ell_{\mathcal{B}_\mathcal{P}^\upsilon}(\boldsymbol{z}).
\end{equation*}
Taking into account that $0 < \tau^{\max}_{\mathcal{B}_\mathcal{P}, \ell}(\vartheta) < \pi - \vartheta/2$, we transform $\tau \in (0, \tau^{\max}_{\mathcal{B}_\mathcal{P}, \ell}(\vartheta))$ into
\begin{equation*}
    -2\cos\frac\vartheta2 < -2\cos\left(\tau + \frac\vartheta2\right) < -2\cos\left(\tau^{\max}_{\mathcal{B}_\mathcal{P}, \ell}(\vartheta) + \frac\vartheta2\right).
\end{equation*}
Substituting $\vartheta$ in the first equation and using obtained inequality we conclude that
\begin{multline*}
    \mathcal{F}_{\mathcal{B}_\mathcal{P}^\upsilon}^\ell \overset{\mathrm{def}}{=} \left\{\boldsymbol{z} \in \mathbb{R}^2 \times \mathbb{S}:\: 0 < -x\cos\frac\theta2 + y\sin\frac\theta2 \right.\\
    \left. < 2\upsilon\left(\cos\frac\theta2 -\cos\left(\tau^{\max}_{\mathcal{B}_\mathcal{P}^\upsilon, \ell}(\theta) + \frac{\upsilon\theta}2\right)\right), \: 0 < \theta < 2\pi\right\},
\end{multline*}
where
\begin{equation*}
    \tau^{\max}_{\mathcal{B}_\mathcal{P}^\upsilon, \ell}(\theta) \overset{\mathrm{def}}{=} \tau^{\max}_{\mathcal{B}_\mathcal{P}, \ell}\left((1 - \upsilon)\pi + \upsilon\theta\right).
\end{equation*}

To obtain the same parametric description of $\mathcal{B}_\mathcal{TS}^\upsilon$ we must solve the system 
\begin{equation*}
    \left\{
    \begin{aligned}
        &\bar{\boldsymbol{z}}_{\mathcal{B}_\mathcal{TS}^\upsilon}(\tau, \vartheta) = \boldsymbol{z};\\
        &\vartheta \in (0, \vartheta^{12}_\ell), \: \tau \in \left(\vartheta, \tau_{\mathcal{B}_\mathcal{TS}, \ell}^{\max}(\vartheta)\right).
    \end{aligned}
    \right.
\end{equation*}
Expressing the $\sin(\tau - \vartheta)$ and $\cos(\tau - \vartheta)$ from the system gives 
\begin{equation*}
    \sin(\tau - \vartheta) = \frac{\upsilon x - 1 + \cos\vartheta}{\ell + \vartheta}, \quad \cos(\tau - \vartheta) = \frac{y + \sin\vartheta}{\ell + \vartheta}.
\end{equation*}
For a fixed $\vartheta \in (0, \vartheta^{12}_\ell)$, it corresponds to a parametric description of a circle on $xy$-plane where $\tau$ is a parameter. Using $\tau > \vartheta$ and Lemma~\ref{lem:tau_max_BTD} we obtain
\begin{equation*}
    \sin(\tau - \vartheta) > 0, \quad \cos(\tau - \vartheta) > \cos(\tau_{\mathcal{B}_\mathcal{TS}, \ell}^{\max}(\vartheta) - \vartheta).
\end{equation*}
Expressing $\vartheta$ and substituting in these inequalities we obtain the frame
\begin{multline*}
    \mathcal{F}_{\mathcal{B}_\mathcal{TS}^\upsilon}^\ell \overset{\mathrm{def}}{=} \{\boldsymbol{z} \in \mathbb{R}^2 \times \mathbb{S}:\: 0 < (1 + \upsilon)\pi - \upsilon\theta < \vartheta^{12}_\ell, \: 1 - \cos\theta < \upsilon x,\\
    (\ell + (1 + \upsilon)\pi - \upsilon\theta)\cos(\tau^{\max}_{\mathcal{B}_\mathcal{TS}^\upsilon, \ell}(\theta) + \upsilon\theta) < y - \upsilon\sin\theta\},
\end{multline*}
where
\begin{equation*}
    \tau^{\max}_{\mathcal{B}_\mathcal{TS}^\upsilon, \ell}(\theta) \overset{\mathrm{def}}{=} \tau^{\max}_{\mathcal{B}_\mathcal{TS}, \ell}\left((1 + \upsilon)\pi - \upsilon\theta\right).
\end{equation*}
Eliminating $\tau$ from the circle parametrization and expressing $\ell$ we have
\begin{equation*}
    \ell = -(1 + \upsilon)\pi + \upsilon\theta + \sqrt{\left(\upsilon x - 1 + \cos\theta\right)^2 + (y - \upsilon\sin\theta)^2} \overset{\mathrm{def}}{=} \ell_{\mathcal{B}_\mathcal{TS}^\upsilon}(\boldsymbol{z}).
\end{equation*}

To obtain the parametric description of $\mathcal{B}_\mathcal{TD}^{\upsilon}$ we have to exclude $\tau$ and $\vartheta$ from the system of equations
\begin{equation*}
    \left\{
    \begin{aligned}
        &\bar{\boldsymbol{z}}_{\mathcal{B}_\mathcal{TD}^\upsilon}(\tau, \vartheta) = \boldsymbol{z};\\
        &\vartheta \in (0, 2\pi - \vartheta^{21}_\ell), \quad \tau \in \left(\vartheta/2, \tau_{\mathcal{B}_\mathcal{TD}, \ell}^{\max}(\vartheta)\right).
    \end{aligned}
    \right.
\end{equation*}
Expressing the $\sin(\tau - \vartheta)$ and $\cos(\tau - \vartheta)$ from the system gives
\begin{align*}
    \sin(\tau - \vartheta) = \frac{(\upsilon x + 1 + \cos\vartheta)(\ell + 2\tau - \vartheta) - 2(y + \sin\vartheta)}{(\ell + 2\tau - \vartheta)^2 + 4},\\
    \cos(\tau - \vartheta) = \frac{2(\upsilon x + 1 + \cos\vartheta) + (y + \sin\vartheta)(\ell + 2\tau - \vartheta)}{(\ell + 2\tau - \vartheta)^2 + 4}.
\end{align*}
The sum of squares of these values and inequality $\ell + 2\tau - \vartheta > 0$ give
\begin{equation}\label{eq:ell2tauvartheta}
    \ell + 2\tau - \vartheta = \sqrt{(\upsilon x + 1 + \cos\vartheta)^2 + (y + \sin\vartheta)^2 - 4}.
\end{equation}
Since $\tau - \vartheta \in (-\vartheta/2, \tau_{\mathcal{B}_\mathcal{TD}, \ell}^{\max}(\vartheta) - \vartheta) \subset (-\pi, 0)$, $\sin(\tau - \vartheta)$ must be negative. Hence, we can inverse $\cos(\tau - \vartheta)$ as
\begin{equation}\label{eq:tauvartheta}
    \tau - \vartheta = -\arccos\frac{2(\upsilon x + 1 + \cos\vartheta) + (y + \sin\vartheta)(\ell + 2\tau - \vartheta)}{(\ell + 2\tau - \vartheta)^2 + 4}.
\end{equation}
Combining \eqref{eq:ell2tauvartheta} and \eqref{eq:tauvartheta} and using $\vartheta = (1 + \upsilon)\pi - \upsilon\theta$ we obtain 
\begin{multline*}
    \ell = \sqrt{(\upsilon x + 1 + \cos\theta)^2 + (y - \upsilon\sin\theta)^2 - 4} - (1 + \upsilon)\pi + \upsilon\theta\\
    + 2\arccos\left(\frac{2(\upsilon x + 1 + \cos\theta)}{(\upsilon x + 1 + \cos\theta)^2 + (y -\upsilon\sin\theta)^2}\right.\\
    \left. + \frac{(y - \upsilon\sin\theta)\sqrt{(\upsilon x + 1 + \cos\theta)^2 + (y - \upsilon\sin\theta)^2 - 4}}{(\upsilon x + 1 + \cos\theta)^2 + (y -\upsilon\sin\theta)^2}\right) \overset{\mathrm{def}}{=} \ell_{\mathcal{B}_\mathcal{TD}^\upsilon}(\boldsymbol{z}).
\end{multline*}
Given the above, the frame is described by
\begin{multline*}
    \mathcal{F}_{\mathcal{B}_\mathcal{TD}^\upsilon}^\ell \overset{\mathrm{def}}{=} \{ \boldsymbol{z} \in \mathbb{R}^2 \times \mathbb{S}:\: 0 < (1 + \upsilon)\pi - \upsilon\theta < 2\pi - \vartheta^{21}_\ell,\\
    0 < -\ell + \sqrt{(\upsilon x + 1 + \cos\theta)^2 + (y - \upsilon\sin\theta)^2 - 4}\\
    < 2\tau^{\max}_{\mathcal{B}_\mathcal{TD}^\upsilon, \ell}(\theta) - (1 + \upsilon)\pi + \upsilon\theta,\\
    (\upsilon x + 1 + \cos\theta)\sqrt{(\upsilon x + 1 + \cos\theta)^2 + (y - \upsilon\sin\theta)^2 - 4} < 2(y - \upsilon\sin\theta)\},
\end{multline*}
where
\begin{equation*}
    \tau^{\max}_{\mathcal{B}_\mathcal{TD}^\upsilon, \ell}(\theta) \overset{\mathrm{def}}{=} \tau^{\max}_{\mathcal{B}_\mathcal{TD}, \ell}((1 + \upsilon)\pi - \upsilon\theta).
\end{equation*}

The parametric descriptions \eqref{eq:BPL_BUL_vartheta} of $\mathcal{B}_\mathcal{PL}^\upsilon$ and $\mathcal{B}_\mathcal{UL}^\upsilon$ contain only one parameter $\vartheta$ which can be expressed through $\theta$ in both cases. Hence,
\begin{align*}
    &\mathcal{B}_\mathcal{PL}^\upsilon = \{\boldsymbol{z} \in \mathbb{R}^2 \times \mathbb{S}:\: \bar{\boldsymbol{z}}_{\mathcal{B}_\mathcal{PL}^\upsilon}((1 - \upsilon)\pi + \upsilon\theta) = \boldsymbol{z}, \: \vartheta^{21}_\ell < (1 - \upsilon)\pi + \upsilon\theta < 2\pi\},\\
    &\mathcal{B}_\mathcal{UL}^\upsilon = \{\boldsymbol{z} \in \mathbb{R}^2 \times \mathbb{S}:\: \bar{\boldsymbol{z}}_{\mathcal{B}_\mathcal{UL}^\upsilon}((1 + \upsilon)\pi - \upsilon\theta) = \boldsymbol{z}, \: 0 < (1 + \upsilon)\pi - \upsilon\theta < \vartheta^{12}_\ell\}.
\end{align*}

\begin{proposition}\label{prop:optimal_controls_on_barrier}
    The optimal feedback controls on the barrier of the pursuer and evader in the GTIC are given by
    \begin{align*}
        &u^*(\boldsymbol{z}) =
        \begin{cases}
            0,& \quad \boldsymbol{z} \in \mathcal{B}_\mathcal{UL}^{-1} \cup \mathcal{B}_\mathcal{UL}^{+1};\\
            +1,& \quad \boldsymbol{z} \in \mathcal{B}_\mathcal{P}^{-1} \cup \mathcal{B}_\mathcal{PL}^{-1} \cup \mathcal{B}_\mathcal{TS}^{+1} \cup \mathcal{B}_\mathcal{TD}^{-1} \cup \mathcal{BUP}_{-1} \cup \mathcal{B}_\mathcal{DL};\\
            -1,& \quad \boldsymbol{z} \in \mathcal{B}_\mathcal{P}^{+1} \cup \mathcal{B}_\mathcal{PL}^{+1} \cup \mathcal{B}_\mathcal{TS}^{-1} \cup \mathcal{B}_\mathcal{TD}^{+1} \cup \mathcal{BUP}_{+1} \cup \mathcal{B}_\mathcal{DL};\\
            \mathrm{sgn}\:x,& \quad \boldsymbol{z} \in \mathcal{BUP}_0,
        \end{cases}\\
        &v^*(\boldsymbol{z}) =
        \begin{cases}
            +1,& \quad \boldsymbol{z} \in \mathcal{B}_\mathcal{P}^{+1} \cup \mathcal{B}_\mathcal{PL}^{+1} \cup \mathcal{B}_\mathcal{TS}^{+1} \cup \mathcal{B}_\mathcal{TD}^{+1} \cup \mathcal{B}_\mathcal{UL}^{+1} \cup \mathcal{BUP}_{+1} \cup \mathcal{B}_\mathcal{DL};\\
            -1,& \quad \boldsymbol{z} \in \mathcal{B}_\mathcal{P}^{-1} \cup \mathcal{B}_\mathcal{PL}^{-1} \cup \mathcal{B}_\mathcal{TS}^{-1} \cup \mathcal{B}_\mathcal{TD}^{-1} \cup \mathcal{B}_\mathcal{UL}^{-1} \cup \mathcal{BUP}_{-1} \cup \mathcal{B}_\mathcal{DL};\\
            \mathrm{sgn}\:x,& \quad \boldsymbol{z} \in \mathcal{BUP}_0.
        \end{cases}
    \end{align*}
\end{proposition}

Proposition \ref{prop:optimal_controls_on_barrier} gives the optimal feedback controls on the barrier. Since the barrier is a two-dimensional manifold in three-dimensional space, the direct calculations with floating points will almost always miss the barrier. This problem can be avoided by using the special representation of all pieces of the barrier \eqref{eq:special_form_of_pieces}. The natural way is to check to belong to the layer
\begin{equation*}
    \mathcal{P}' = \{\boldsymbol{z} \in \mathcal{F}^{\ell_\mathcal{P}(\boldsymbol{z})}_\mathcal{P}:\: \ell \leq \ell_\mathcal{P}(\boldsymbol{z}) \leq \ell(1 + \delta)\}
\end{equation*}
for each piece of the barrier. Here $\delta \in \mathbb{R}^+$ is a relative layer width. If calculations are precise, then the evader can guarantee the evasion with minimal approaching distance $\ell(1 + \delta)$.

\section{Conclusions}\label{sec:conclusions}

In this paper, we derived a complete analytical description of the barrier for the GTIC. The obtained description is provided for all possible values of the capture radius. An analysis of the variation in the barrier geometry with a variation in the capture radius showed that there are differences between the shape of the barrier for "small" and "large" values of the capture radius. These differences lead to the fact that the analytical conditions for cutting off redundant parts of semipermeable surfaces that constitute the barrier have different forms.

Taking into account the classification of possible forms of the barrier depending on the capture radius, it is no longer necessary to use the trial and error method to construct the barrier, as was done in other approaches when analyzing the barrier of the GTIC. In addition, in this paper, we obtained a parametric description of the barrier in terms of the state vector. Such a description makes it possible to synthesize feedback optimal controls that do not require the elimination of additional parameters by numerical methods, so the resulting controls can be calculated directly from the state vector.

For the collision avoidance problem, the obtained expression for the evader's optimal control on the barrier can be made resistant to the rounding error of the state vector. To do this, instead of calculating the control for the barrier surface, we should calculate the control by belonging to the layer whose boundary is the barrier. Belonging to such a layer can also be established by analytical calculations. This circumstance makes the evasive control more attractive for practical implementation.

One possible avenue for further research is to investigate the barrier of the GTC by the same methods.

\bibliographystyle{spmpsci_unsrt}
\bibliography{main.bib}

\end{document}